\documentclass[11pt,letterpaper,oneside]{article}

\usepackage{latexsym}
\usepackage{amsmath}
\usepackage{amsthm}
\usepackage{amssymb}
\usepackage{amsfonts}
\usepackage{fancyhdr}
\usepackage{comment}

\usepackage{mathrsfs}                           % For different calligraphic fonts

\newcommand{\mBic}{\,\Leftrightarrow\,}         % Formatting for the biconditional

\newcommand{\mR}{\mathbf{R}}                    % Formatting for R
\newcommand{\mC}{\mathbf{C}}                    % Formatting for C
\newcommand{\abs}[1]{\lvert #1 \rvert}          % Formatting for the absolute value
\newcommand{\norm}[1]{\lVert #1 \rVert}         % Formatting for the norm
\newcommand{\br}[1]{\langle #1 \rangle}         % Formatting for the inner product
\newcommand{\mS}{\mathscr{S}}

\newcommand{\mSp}{\mathscr{S}^{\prime}}

\newcommand{\mF}{\mathscr{F}}

\newcommand{\ehat}{\,\hat{\rule{0pt}{6pt}}\,}

\newcommand{\re}{\mathrm{Re}}
\newcommand{\im}{\mathrm{Im}}

\newcommand{\curl}{\mathrm{curl}}

\newcommand{\closure}[1]{\overline{#1}}
\newcommand{\dbar}{\overline{\partial}}

\newcommand{\mOp}{\mathrm{Op}}

\theoremstyle{definition}
\newtheorem{thm}{Theorem}[section]
\newtheorem{prop}{Proposition}[section]

\newtheorem{lemma}{Lemma}[section]
\newtheorem*{definition}{Definition}

\newtheorem*{remark}{Remark}

\title{Semiclassical Pseudodifferential Calculus and the \\
Reconstruction of a Magnetic Field}
\author{ \\
Mikko Salo \\[5pt]
{\small Department of Mathematics and Statistics / RNI} \\
{\small P.O. Box 68, 00014 University of Helsinki, Finland} \\
{\small mikko.salo@helsinki.fi}}
\date{}

\begin{document}

\maketitle

\begin{abstract}
We give a procedure for reconstructing a magnetic field and electric potential from boundary measurements given by the Dirichlet to Neumann map for the magnetic Schr\"odinger operator in $\mR^n$, $n \geq 3$. The magnetic potential is assumed to be continuous with $L^{\infty}$ divergence and zero boundary values. The method is based on semiclassical pseudodifferential calculus and the construction of complex geometrical optics solutions in weighted Sobolev spaces.
\end{abstract}

\section{Introduction}

Let $\Omega \subseteq \mR^n$, $n \geq 3$, be a bounded domain with $C^{1,1}$ boundary. We consider the magnetic Schr\"odinger operator 
\begin{equation*}
H_{W,q} = \sum_{j=1}^n (D_j + W_j)^2 + q
\end{equation*}
where $D_j = \frac{1}{i} \frac{\partial}{\partial x_j}$, $W \in L^{\infty}(\Omega ; \mC^n)$ is the magnetic potential, and $q \in L^{\infty}(\Omega ; \mC)$ is the electric potential (the coefficients can be complex valued). Assuming $0$ is not a Dirichlet eigenvalue of $H_{W,q}$ in $\Omega$, the problem 
\begin{equation*} \label{auxiliarydirichlet}
\left\{ \begin{array}{rll}
H_{W,q} u &\!\!\!= 0 & \quad \text{in } \Omega, \\
u &\!\!\!= f & \quad \text{on } \partial \Omega
\end{array} \right.
\end{equation*}
has a unique solution $u = u_f \in H^1(\Omega)$ for any $f \in H^{1/2}(\partial \Omega)$.

The boundary measurements are given by the Dirichlet to Neumann map (DN map), defined formally by 
\begin{equation*}
\Lambda_{W,q}: f \mapsto\frac{\partial u_f}{\partial \nu} \Big|_{\partial \Omega} + i (W \cdot \nu) f.
\end{equation*}
Here $\nu$ is the outer unit normal to $\partial \Omega$. More precisely, if $f, g \in H^{1/2}(\partial \Omega)$ we define $\Lambda_{W,q}$ using the equivalent weak formulations 
\begin{align*}
\langle \Lambda_{W,q} f, g \rangle &= \int_{\Omega} (\nabla u_f \cdot \nabla e_g + iW \cdot (u_f \nabla e_g - e_g \nabla u_f) + (W^2 + q) u_f e_g) \,dx \\
 &= \int_{\Omega} (\nabla e_f \cdot \nabla v_g + iW \cdot (e_f \nabla v_g - v_g \nabla e_f) + (W^2 + q) e_f v_g) \,dx
\end{align*}
where $v_h \in H^1(\Omega)$ solves the adjoint problem $H_{-W,q} v_h = 0$ in $\Omega$ with $v_h|_{\partial \Omega} = h$, and $e_h$ is any $H^1(\Omega)$ function with $e_h|_{\partial \Omega} = h$. Then $\Lambda_{W,q}$ is a bounded map $H^{1/2}(\partial \Omega) \to H^{-1/2}(\partial \Omega)$.

The gauge transformation $W \mapsto W + \nabla p$, where $p \in W^{1,\infty}(\Omega ; \mC)$, transforms the magnetic potential to a gauge equivalent potential but preserves the magnetic field $\curl\,W$. If additionally $p|_{\partial \Omega} = 0$ then $\Lambda_{W+\nabla p,q} = \Lambda_{W,q}$, which means that boundary measurements are preserved in gauge transformations which respect the boundary.

We are interested in recovering $\curl\,W$ and $q$ from $\Lambda_{W,q}$. This is a typical inverse problem where one wishes to know the interior properties of a medium by making measurements at the boundary. It is related to the extensively studied inverse conductivity problem of Calder{\'o}n \cite{calderon}. In fact, most known results in dimensions $n \geq 3$ (starting in Sylvester-Uhlmann \cite{sylvesteruhlmann}, for later work see \cite{uhlmannasterisque}, \cite{uhlmanndevelopments}) reduce that problem to recovering $q$ from $\Lambda_{0,q}$.

The problem is closely related to inverse scattering at fixed energy, which was studied for the magnetic case in Eskin-Ralston \cite{eskinralston}. See also Sun \cite{sunnote} for an application of the method in \cite{eskinralston}. One motivation for the present study has been to understand the approach of \cite{eskinralston} and to clarify the role of pseudodifferential operators in inverse problems for first order perturbations of the Laplacian.

Previous results for this inverse problem concern unique determination of the coefficients, and they state that $\Lambda_{W_1,q_1} = \Lambda_{W_2,q_2}$ implies $\curl\,W_1 = \curl\,W_2$ and $q_1 = q_2$ in $\Omega$, under varying assumptions on $W_j$, $q_j$, and $\Omega$. Sun \cite{sun} proved this in the case where the $\curl\,W_j$ are small. Panchenko \cite{panchenko} proves a similar result in a less regular setting. The smallness assumption was removed by Nakamura-Sun-Uhlmann \cite{nakamurasunuhlmann} who considered $C^{\infty}$ coefficients. Based on the method in \cite{nakamurasunuhlmann}, the smoothness assumption was reduced to $C^1$ by Tolmasky \cite{tolmasky}, and to Dini continuous in \cite{salo}. All these results are nonconstructive. We also mention the recent papers \cite{brownsalo}, \cite{dksu} and \cite{tzou} which consider boundary determination, partial Cauchy data and stability for this inverse problem.

In this paper we give a constructive algorithm for recovering $\curl\,W$ and $q$ from $\Lambda_{W,q}$. Such an algorithm is well known in the case $W = 0$ and is due to Nachman \cite{nachman} and also Novikov \cite{novikov}. This result combines ideas from scattering theory with so called complex geometrical optics (CGO) solutions of the Schr\"odinger equation, which were introduced in the fundamental paper \cite{sylvesteruhlmann}. The main point is that the CGO solutions are defined globally and are unique in a certain sense. The problem in extending the results of Nachman to nonzero $W$ has been that the main method of producing CGO solutions in that case, the pseudodifferential conjugation technique of Nakamura-Uhlmann \cite{nakamurauhlmann}, only works in bounded domains and there does not seem to be a proper notion of uniqueness for solutions.

We will give a global version of the Nakamura-Uhlmann technique which will produce global CGO solutions with uniqueness in the proper weighted Sobolev spaces. To do this we apply semiclassical pseudodifferential calculus. This is largely equivalent with the parameter-dependent calculus used earlier in such results, but it simplifies the proofs. A main new element in our approach is a variant of the pseudodifferential cutoff technique used by Takeuchi \cite{takeuchi} and Kenig-Ponce-Vega \cite{kenigponcevega} in the context of nonlinear Schr\"odinger equations. We remark that semiclassical notation was also recently used by Kenig-Sj\"ostrand-Uhlmann \cite{kenigsjostranduhlmann} who studied the problem of recovering coefficients from partial boundary measurements.

We record some notation. Let $\Delta = \sum D_j^2$, $\Delta_{\zeta} = \Delta + 2\zeta \cdot D$, and $D_{\zeta} = D + \zeta$, where $D = (D_1,\ldots,D_n)$ and $\zeta \in \mC^n$, $\zeta^2 := \zeta \cdot \zeta = 0$. The $\zeta$-dependent operators $\Delta_{\zeta}$ and $D_{\zeta}$ arise naturally in the construction of CGO solutions. For $\delta \in \mR$ we also use weighted $L^2$ spaces $L^2_{\delta}$ with norm $\norm{f}_{L^2_{\delta}} = \norm{\br{x}^{\delta} f}_{L^2}$ where $\br{x} = (1+\abs{x}^2)^{1/2}$, and weighted Sobolev spaces $H^s_{\delta}$ with norm $\norm{f}_{H^s_{\delta}} = \norm{\br{x}^{\delta} f}_{H^s}$. If $X$ is a function space we write $X_c$ for the set of compactly supported functions in $X$, and $X_{\Omega}$ for all functions in $X$ supported in $\closure{\Omega}$.

The construction of CGO solutions follows from the norm estimates in the following theorem.

\begin{thm} \label{normestimates}
Let $W \in C_c(\mR^n;\mC^n)$, $q \in L^{\infty}_c(\mR^n ; \mC)$, and $-1 < \delta < 0$. If $\zeta \in \mC^n$ with $\zeta^2 = 0$ and $\abs{\zeta}$ is large enough, then for any $f \in L^2_{\delta+1}(\mR^n)$ the equation 
\begin{equation} \label{mainequation}
(\Delta_{\zeta} + 2 W \cdot D_{\zeta} + q)u = f
\end{equation}
has a unique solution $u \in H^1_{\delta}(\mR^n)$. Furthermore, $u \in H^2_{\delta}(\mR^n)$, and $u$ satisfies for $0 \leq s \leq 2$ 
\begin{equation*}
\norm{u}_{H^s_{\delta}} \leq C \abs{\zeta}^{s-1} \norm{f}_{L^2_{\delta+1}},
\end{equation*}
where $C$ is independent of $\zeta$ and $f$.
\end{thm}

Theorem \ref{normestimates} generalizes results in \cite{sylvesteruhlmann} to operators with first order terms. Using the global CGO solutions to $H_{W,q} u = 0$ obtained from this theorem, we may extend the results of Nachman to obtain a constructive algorithm for recovering $\curl\,W$ and $q$ from $\Lambda_{W,q}$. The main result is as follows.

\begin{thm} \label{reconstruction}
Let $\Omega \subseteq \mR^n$, $n \geq 3$, be a bounded simply connected $C^{1,1}$ domain. Suppose $W \in C_{\Omega}(\mR^n ; \mC^n)$ with $D \cdot W \in L^{\infty}(\mR^n ; \mC)$, and suppose $q \in L^{\infty}_{\Omega}(\mR^n ; \mC)$. Also suppose that $0$ is not a Dirichlet eigenvalue of $H_{W,q}$ in $\Omega$. Then $\Lambda_{W,q}$ determines $\curl\,W$ uniquely and constructively. Further, if $W$ is $C^{1+\varepsilon}$ and $\partial \Omega$ is $C^{2+\varepsilon}$ for some $\varepsilon > 0$, then one may construct $q$ from $\Lambda_{W,q}$.
\end{thm}

The reconstruction procedure for $\curl\,W$ is outlined in the following four steps.

\begin{enumerate}
\item 
From the knowledge of $\Lambda_{W,q}$, one may determine the boundary values $u_{\zeta}|_{\partial \Omega}$ of a CGO solution $u_{\zeta}$ as the unique solution of a boundary integral equation on $\partial \Omega$.
\item 
From $\Lambda_{W,q}$ and $u_{\zeta}|_{\partial \Omega}$ one computes a scattering transform $t_{W,q}(\xi,\zeta)$.
\item 
The expression $R_{W,q}(\xi,\mu) = \lim_{s\to\infty} s^{-1} t_{W,q}(\xi,s\mu)$ is essentially the Fourier transform of $\curl\,W$.
\item 
$\curl\,W$ may be computed from $R_{W,q}$ using the inverse Fourier transform.
\end{enumerate}

The structure of the paper is as follows. Section 2 contains some facts on semiclassical pseudodifferential calculus, and Section 3 contains estimates for $\dbar$ equations which will be needed later. In Section 4 we prove Theorem \ref{normestimates}, where the main step is to conjugate the first order term into a lower order one using pseudodifferential operators. Section 5 discusses equivalent problems  which characterize the CGO solutions. In Sections 6 and 7 we reconstruct the magnetic field and electric potential, respectively.

\section{Semiclassical pseudodifferential calculus}

Our method will involve pseudodifferential operators depending on a small parameter $h$. We will collect the required properties of these operators here. See \cite{dimassisjostrand} for more details.

\begin{definition}
If $0 \leq \sigma < 1/2$ and $m \in \mR$, we let $S^m_{\sigma}(\mR^n)$ be the space of all functions $a(x,\xi;h)$ where $x,\xi \in \mR^n$ and $h \in (0,h_0]$, $h_0 \leq 1$, such that $a(\,\cdot\,;h) \in C^{\infty}(\mR^{2n})$ and 
\begin{equation*}
\abs{\partial_x^{\alpha} \partial_{\xi}^{\beta} a(x,\xi;h)} \leq C_{\alpha \beta} h^{-\sigma\abs{\alpha+\beta}} \br{\xi}^m
\end{equation*}
for all $\alpha,\beta$. If $a \in S^m_{\sigma}$ we define an operator $A = \mOp_h(a) = a(x,hD)$ by 
\begin{equation*}
Af(x) = (2\pi)^{-n} \int_{\mR^n} e^{ix\cdot\xi} a(x,h\xi;h) \hat{f}(\xi) \,d\xi.
\end{equation*}
The class of such operators is denoted by $\mOp\,S^m_{\sigma}$.
\end{definition}

Note that we use the standard quantization instead of Weyl quantization in the definition of the operators.

\begin{prop} \label{psdoproperties}
\cite{dimassisjostrand} Let $a \in S^m_{\sigma}$ with $m \in \mR$ and $0 \leq \sigma < 1/2$.
\begin{enumerate}
\item[(a)] 
$\mOp_h(a)$ is a continuous map $\mS \to \mS$ and $\mSp \to \mSp$.
\item[(b)] 
If $m = 0$ then $\mOp_h(a)$ is bounded $L^2 \to L^2$, and there is a constant $C$ with 
\begin{equation*}
\norm{\mOp_h(a)}_{L^2 \to L^2} \leq C
\end{equation*}
for $0 < h \leq h_0$.
\item[(c)] 
$\partial_{x_j} \mOp_h(a) = \mOp_h(a) \partial_{x_j} + \mOp_h\Big(\frac{\partial a}{\partial x_j}\Big)$.
\item[(d)]
The adjoint $\mOp_h(a)^* = \mOp_h(a^*)$, where $a^* \in S^m_{\sigma}$ satisfies for any $N$ 
\begin{equation*}
a^* = \sum_{\abs{\alpha} < N} \frac{h^{\abs{\alpha}} \partial_{\xi}^{\alpha} D_x^{\alpha} a}{\alpha!} + h^{N(1-2\sigma)} S^m_{\sigma}.
\end{equation*}
\item[(e)] 
If $a \in S^m_{\sigma}$ and $b \in S^{m'}_{\sigma}$ then $\mOp_h(a) \mOp_h(b) = \mOp_h(c)$ where $c \in S^{m+m'}_{\sigma}$ satisfies for any $N$ 
\begin{equation*}
c = \sum_{\abs{\alpha} < N} \frac{h^{\abs{\alpha}} \partial_{\xi}^{\alpha} a D_x^{\alpha} b}{\alpha!} + h^{N(1-2\sigma)} S^{m+m'}_{\sigma}.
\end{equation*}
Also, $[\mOp_h(a), \mOp_h(b)] = \mOp_h(d)$ where $d \in S^{m+m'}_{\sigma}$ and 
\begin{equation*}
d = \frac{h}{i} H_{a} b + h^{2(1-2\sigma)} S^{m+m'}_{\sigma}
\end{equation*}
where $H_a = \nabla_{\xi} a \cdot \nabla_x - \nabla_x a \cdot \nabla_{\xi}$ is the Hamilton vector field of $a$.
\end{enumerate}
\end{prop}

We will also need the boundedness of $\mOp\,S^0_{\sigma}$ operators on weighted Sobolev spaces. Here it is natural to use semiclassical spaces $H^s_{\delta,h}$, defined by $\norm{f}_{H^s_{\delta,h}} = \norm{\br{hD}^s \br{x}^{\delta} f}_{L^2}$. If $s \geq 0$ is an integer an equivalent norm is given by $\sum_{\abs{\alpha} \leq s} h^{\abs{\alpha}} \norm{D^{\alpha} f}_{L^2_{\delta}}$.

\begin{prop} \label{weightedpsdoestimates}
Let $a \in S^0_{\sigma}$ with $0 \leq \sigma < 1/2$ and $s_0, \delta_0 \geq 0$. Then $\mOp_h(a)$ is bounded $H^s_{\delta,h}(\mR^n) \to H^s_{\delta,h}(\mR^n)$ for any $s, \delta \in \mR$, and there is a constant $C$ with $\norm{\mOp_h(a)}_{H^s_{\delta,h} \to H^s_{\delta,h}} \leq C$ whenever $\abs{s} \leq s_0$, $\abs{\delta} \leq \delta_0$, and $0 < h \leq h_0$.
\end{prop}
\begin{proof}
All constants below are independent of $h$. We begin by showing that 
\begin{equation} \label{psdoweightedfirstestimate}
\norm{\br{x}^{\delta} \br{hD}^s Af}_{L^2} \leq C \norm{\br{x}^{\delta} \br{hD}^s f}_{L^2}.
\end{equation}
It is enough to take $s = 0$ since otherwise we may consider $\br{hD}^s A \br{hD}^{-s}$, which is in $\mOp\,S^0_{\sigma}$ by Proposition \ref{psdoproperties}. If $m$ is a nonnegative integer we define $Tf(x) = \br{x}^{-2m} A (\br{x}^{2m} f)$. For $f \in \mS$ one has  
\begin{align*}
Tf(x) &= (2\pi)^{-n} \int e^{ix \cdot \xi} a(x,h\xi) \br{x}^{-2m} (I-\Delta_{\xi})^m \hat{f}(\xi) \,d\xi \\
 &= (2\pi)^{-n} \int (I-\Delta_{\xi})^m ( \br{x}^{-2m} e^{ix \cdot \xi} a(x,h\xi) ) \hat{f}(\xi) \,d\xi.
\end{align*}
It follows by differentiation that $T$ is in $\mOp\,S^0_{\sigma}$, and Proposition \ref{psdoproperties} gives \eqref{psdoweightedfirstestimate} for $\delta = -2m$. The estimate for $\delta \leq 0$ follows from the Stein-Weiss interpolation theorem, and for $\delta \geq 0$ by duality using the fact that $\mOp_h(a)^*$ is a pseudodifferential operator.

It remains to relate \eqref{psdoweightedfirstestimate} to the norm $\norm{\br{hD}^s \br{x}^{\delta} f}_{L^2}$. If $k \geq 0$ is an integer then $\norm{\br{hD}^{2k} \br{x}^{\delta} f}_{L^2} \leq C \sum_{\abs{\alpha} \leq 2k} \norm{\br{x}^{\delta} (hD)^{\alpha} f}_{L^2}$, using $h \leq 1$. We claim that 
\begin{equation} \label{psdoweightedsecondestimate}
\norm{\br{x}^{\delta} (hD)^{\alpha} f}_{L^2} \leq C \norm{\br{x}^{\delta} \br{hD}^{2k} f}_{L^2}, \quad \abs{\alpha} \leq 2k.
\end{equation}
In fact, $\norm{\br{x}^{2m} (hD)^{\alpha} f}_{L^2} = C \norm{\br{D_{\xi}}^{2m} r(\xi) \br{h\xi}^{2k} \hat{f}}_{L^2}$ if $m \geq 0$ is an integer, where $r(\xi)$ has bounded derivatives of all orders. Differentiation gives that this is bounded by $C \sum_{\abs{\beta} \leq 2m} \norm{D_{\xi}^{\beta} \br{h\xi}^{2k} \hat{f}}_{L^2}$. Going back to the $x$-side gives \eqref{psdoweightedsecondestimate} for $\delta = 2m$, and the estimate follows for $\delta \in \mR$ by interpolation and duality. Using \eqref{psdoweightedsecondestimate} implies that 
\begin{equation*}
\norm{\br{hD}^{2k} \br{x}^{\delta} f}_{L^2} \leq C \norm{\br{x}^{\delta} \br{hD}^{2k} f}_{L^2}.
\end{equation*}
The last estimate applied to $Af$ and then \eqref{psdoweightedfirstestimate} give 
\begin{align*}
\norm{\br{hD}^{2k} \br{x}^{\delta} Af}_{L^2} &\leq C \norm{\br{x}^{\delta} \br{hD}^{2k} f}_{L^2} \\
 &\leq C \sum_{\abs{\alpha} \leq 2k} \norm{\br{x}^{\delta} (hD)^{\alpha} f}_{L^2} \\
 &\leq C \sum_{\abs{\alpha} \leq 2k} \norm{(hD)^{\alpha} \br{x}^{\delta} f}_{L^2} \\
 &\leq C \norm{\br{hD}^{2k} \br{x}^{\delta} f}_{L^2}.
\end{align*}
The interpolation (now the Stein-Weiss theorem on the Fourier side) and duality give the desired result.
\end{proof}

\section{Estimates for $\dbar$ equations}

In this section we collect some elementary estimates for equations of $\dbar$ type in $\mR^n$. Let $\mu = \gamma_1 + i \gamma_2$ where $\gamma_j \in \mR^n$, $\abs{\gamma_j} = 1$, and $\gamma_1 \cdot \gamma_2 = 0$. The operator $N_{\mu} = \mu \cdot \nabla$ is just $\partial_{x_1} + i \partial_{x_2}$ in different coordinates, so it has an inverse given by 
\begin{equation*}
N_{\mu}^{-1} f(x) = \frac{1}{2\pi} \int_{\mR^2} \frac{1}{y_1 + i y_2} f(x-y_1 \gamma_1 - y_2 \gamma_2) \,dy_1 \,dy_2.
\end{equation*}
This operator satisfies the following.

\begin{lemma} \label{nmuinv1}
Let $f \in W^{k,\infty}(\mR^n)$ with $f = 0$ for $\abs{x} \geq M$. Then $u = N_{\mu}^{-1} f \in W^{k,\infty}(\mR^n)$ solves the equation $N_{\mu} u = f$ in $\mR^n$ and satisfies for $\abs{\alpha} \leq k$ 
\begin{equation} \label{nmuinvestimate}
\abs{\partial^{\alpha} u(x)} \leq C(M) \norm{\partial^{\alpha} f}_{L^{\infty}} \br{x_T}^{-1} \chi_{B(0,M)}(x_{\perp})
\end{equation}
where $x_T$ is the projection of $x$ to the plane $T = \mathrm{span} \{\gamma_1,\gamma_2\}$, $x_{\perp} = x - x_T$, and $\chi_{B(0,M)}$ is the characteristic function of $B(0,M)$.
\end{lemma}
\begin{proof}
Since the statements are rotation invariant we may assume $\gamma_j = e_j$ (the $j$th coordinate vector) and $N_{\mu} = \dbar = \partial_{x_1} + i \partial_{x_2}$. Then $u = N_{\mu}^{-1} f$ is 
\begin{equation} \label{unormalized}
u(x) = \frac{1}{2\pi} \int_{\mR^2} \frac{1}{y_1 + i y_2} f(x'-y',x'') \,dy'.
\end{equation}
We have $u \in L^{\infty}$ since $f \in L^{\infty}_c$. If $\varphi \in C^{\infty}$ satisfies $\varphi = 0$ for $\abs{x} \geq R$ we have 
\begin{equation*}
\langle u, -\dbar \varphi \rangle = -\frac{1}{2\pi} \int_{\mR^n} \int_{\mR^2} \frac{1}{y_1 + i y_2} f(x'-y',x'') \dbar \varphi(x) \,dy' \,dx.
\end{equation*}
The integrand is nonzero only if $\abs{y'} \leq M + R$, which justifies using Fubini's theorem. A change of variables and another use of Fubini's theorem gives 
\begin{equation*}
\langle u, -\dbar \varphi \rangle = \int_{\mR^n} f(x) \Big( -\frac{1}{2\pi} \int_{\mR^2} \frac{1}{y_1 + i y_2} \dbar \varphi(x'+y',x'') \,dy' \Big) \,dx.
\end{equation*}
The inner integral is $\varphi(x)$, which shows that $\dbar u = f$ even when $f \in L^{\infty}_c$.

It is enough to prove \eqref{nmuinvestimate} for $\alpha = 0$. From \eqref{unormalized} we see that $u(x) = 0$ for $\abs{x''} \geq M$. If $\abs{x'} \geq 2M$ then $\abs{y'} \geq \abs{x'}/2$ on the support of the integrand in \eqref{unormalized}, and \eqref{nmuinvestimate} follows. Since \eqref{nmuinvestimate} is easy when $\abs{x'} \leq 2M$ we obtain the required result.
\end{proof}

We will need a version of Lemma \ref{nmuinv1} where $f$ and $\mu$ depend on a para\-meter. Let $V \subseteq \mR^n$ be an open set and let $\gamma_j(\xi)$ ($j=1,2$) be $C^{\infty}$ functions of $\xi \in V$ which satisfy 
\begin{equation} \label{gammaestimates}
1-\varepsilon \leq \abs{\gamma_j(\xi)} \leq 1+\varepsilon, \quad \abs{\gamma_1(\xi) \cdot \gamma_2(\xi)} \leq \varepsilon
\end{equation}
and also $\abs{\partial^{\alpha} \gamma_j(\xi)} \leq M_1$ for $\abs{\alpha} \geq 1$.

\begin{lemma} \label{dbarmodified}
Let $\varepsilon > 0$ be small enough and let $f(x,\xi) \in C^{\infty}(\mR^n \times V)$ satisfy $f(x,\xi) = 0$ for $\abs{x} \geq M$. Then the function 
\begin{equation*}
u(x,\xi) = \frac{1}{2\pi} \int_{\mR^2} \frac{1}{y_1 + i y_2} f(x-y_1 \gamma_1(\xi) - y_2 \gamma_2(\xi), \xi) \,dy_1 \,dy_2
\end{equation*}
is in $C^{\infty}(\mR^n \times V)$, solves $(\gamma_1(\xi) + i\gamma_2(\xi)) \cdot \nabla_x u = f$ in $\mR^n$, and satisfies 
\begin{equation*}
\abs{\partial_x^{\alpha} \partial_{\xi}^{\beta} u(x,\xi)} \leq C_{\alpha \beta M M_1} \Big( \sum_{\abs{\gamma+\delta} \leq \abs{\alpha+\beta}} \norm{\partial_x^{\gamma} \partial_{\xi}^{\delta} f}_{L^{\infty}(\mR^n \times V)} \Big) \br{x_T}^{\abs{\beta}-1} \chi_{B(0,M)}(x_{\perp})
\end{equation*}
where $x_T$ is the projection of $x$ to the plane $T = \mathrm{span} \{\gamma_1(\xi),\gamma_2(\xi)\}$ and $x_{\perp} = x - x_T$.
\end{lemma}
\begin{proof}
Since $f$ is smooth and compactly supported in $x$ it is easy to see that $u$ is smooth and solves the given equation. We have 
\begin{multline*}
\partial_{\xi_j} (f(x-y_1 \gamma_1(\xi) - y_2 \gamma_2(\xi), \xi)) = \nabla_x f(x-y_1 \gamma_1(\xi) - y_2 \gamma_2(\xi), \xi) \cdot \\
 (-y_1 \partial_{\xi_j} \gamma_1(\xi) - y_2 \partial_{\xi_j} \gamma_2(\xi)) + \partial_{\xi_j} f(x-y_1 \gamma_1(\xi) - y_2 \gamma_2(\xi), \xi).
\end{multline*}
Induction and the estimate on the derivatives of $\gamma_j(\xi)$ imply that 
\begin{equation} \label{uest}
\abs{\partial_x^{\alpha} \partial_{\xi}^{\beta} u(x,\xi)} \leq C_{\alpha \beta M_1} \sum_{\abs{\gamma+\delta} \leq \abs{\alpha+\beta}} \sum_{r=0}^{\abs{\beta}} \int_{\mR^2} \abs{y'}^{r-1} \abs{\partial_x^{\gamma} \partial_{\xi}^{\delta} f(x-y_1 \gamma_1(\xi) - y_2 \gamma_2(\xi), \xi)} \,dy'.
\end{equation}
The integrals in \eqref{uest} are over the set 
\begin{equation*}
K(x,\xi) = \{ y' \in \mR^2 \,;\, x-y_1 \gamma_1(\xi) - y_2 \gamma_2(\xi) \in B(0,M) \}.
\end{equation*}
We first note that if $\abs{x_{\perp}} \geq M$ then the right hand side of \eqref{uest} is zero, so we may assume $\abs{x_{\perp}} \leq M$. For the behaviour in $x_T$ we note that 
\begin{align}
y' \in K(x,\xi) &\mBic x - A(y_1,y_2,0)^t \in B(0,M) \notag \\
 &\mBic (y_1,y_2,0)^t \in A^{-1} x + A^{-1} B(0,M), \label{yprimeequiv}
\end{align}
where $A = A(\xi) = (\gamma_1(\xi), \ldots, \gamma_n(\xi))$ is a matrix written in terms of column vectors, and where $\gamma_j(\xi)$ ($3 \leq j \leq n$) are any orthonormal basis of $\{\gamma_1(\xi), \gamma_2(\xi)\}^{\perp}$.

We need to estimate the matrix norm $\norm{A^{-1}} = \sup_{\abs{x}=1} \abs{A^{-1}x}$. From \eqref{gammaestimates} we obtain $\norm{A} \leq 2$ for small $\varepsilon$, so $\abs{A^{-1}x} \geq \abs{x}/2$. An easy calculation using the orthogonality properties of the $\gamma_j$ gives in terms of row vectors 
\begin{equation*}
A^{-1} = (a\gamma_1^t + b\gamma_2^t,\ c\gamma_1^t + d\gamma_2^t,\ \gamma_3^t,\ \ldots,\ \gamma_n^t)
\end{equation*}
where $a,b,c,d$ are obtained from 
\begin{equation*}
\left( \begin{array}{cc} a & b \\ c & d \end{array} \right) = \frac{1}{\abs{\gamma_1}^2 \abs{\gamma_2}^2 - (\gamma_1 \cdot \gamma_2)^2} \left( \begin{array}{cc} \abs{\gamma_2}^2 & -\gamma_1 \cdot \gamma_2 \\ -\gamma_1 \cdot \gamma_2 & \abs{\gamma_1}^2 \end{array} \right).
\end{equation*}
It follows from \eqref{gammaestimates} that for small $\varepsilon$ one has $a, d \sim 1$ and $b, c \sim 0$, implying that $\norm{A^{-1}} \leq 3/2$ for small $\varepsilon$.

Suppose $\abs{x_T} \geq 12M$. If $y' \in K(x,\xi)$ then \eqref{yprimeequiv} gives 
\begin{equation*}
\abs{y'} \geq \abs{A^{-1} x_T} - \abs{A^{-1} x_{\perp}} - \norm{A^{-1}} M \geq \abs{x_T}/2 - 3M/2 - 3M/2 \geq \abs{x_T}/4
\end{equation*}
and also 
\begin{equation*}
\abs{y'} \leq \abs{A^{-1} x_T} + \abs{A^{-1} x_{\perp}} + \norm{A^{-1}} M \leq 3\abs{x_T}/2 + 3M \leq 2\abs{x_T}.
\end{equation*}
Since \eqref{yprimeequiv} implies that $K(x,\xi)$ is contained in a ball of radius $3M/2$, we obtain the desired estimate for $u$ from \eqref{uest} in the case $\abs{x_T} \geq 12M$. If $\abs{x_T} \leq 12M$ then \eqref{yprimeequiv} implies $\abs{y'} \leq 21M$, and the desired estimate follows also in this case.
\end{proof}

\section{Proof of Theorem \ref{normestimates}}

We want to prove the existence, uniqueness, and norm estimates of solutions of 
\begin{equation*}
(\Delta_{\zeta} + 2W \cdot D_{\zeta} + q) u = f
\end{equation*}
with $f \in L^2_{\delta+1}(\mR^n)$. This will be based on the following fundamental estimates for the inverse of $\Delta_{\zeta}$, which imply Theorem \ref{normestimates} in the case $W = q = 0$.

\begin{prop} \label{deltarhoinvestimates}
Let $-1 < \delta < 0$ and let $\zeta \in \mC^n$, $\zeta \cdot \zeta = 0$, $\abs{\zeta} \geq 1$. Then for any $f \in L^2_{\delta+1}(\mR^n)$ the equation 
\begin{equation*}
\Delta_{\zeta} u = f
\end{equation*}
has a unique solution $u \in L^2_{\delta}(\mR^n)$. The solution operator, denoted by $\Delta_{\zeta}^{-1}$, is a bounded map from $L^2_{\delta+1}$ to $H^2_{\delta}$ and satisfies for $0 \leq s \leq 2$ 
\begin{equation*}
\norm{\Delta_{\zeta}^{-1} f}_{H^s_{\delta}} \leq C_0 \abs{\zeta}^{s-1} \norm{f}_{L^2_{\delta+1}},
\end{equation*}
where $C_0 = C_0(n,\delta)$.
\end{prop}
\begin{proof}
The main estimate is the case $s = 0$, which is proved in \cite{sylvesteruhlmann}. Since we could not find a reference for the $H^2_{\delta}$ result we will give the proof here following the $s=1$ case in \cite{brownsingular}.

Take $\phi(\xi) \in C^{\infty}_c(\mR^n)$ with $\phi = 1$ for $\abs{\xi} \leq 4\abs{\zeta}$, $\phi = 0$ for $\abs{\xi} \geq 8\abs{\zeta}$, and $\abs{\nabla \phi} \leq C/\abs{\zeta}$. For $f \in L^2_{\delta+1}$ we write 
\begin{equation} \label{htwodecomposition}
D_j D_k \Delta_{\zeta}^{-1} f = T(\Delta_{\zeta}^{-1} f) + Sf
\end{equation}
where 
\begin{align*}
Tu &= \mF^{-1} \{ \xi_j \xi_k \phi(\xi) \hat{u}(\xi) \}, \\
Sf &= \mF^{-1} \Big\{ \frac{\xi_j \xi_k (1-\phi(\xi))}{\abs{\xi}^2 + 2\zeta \cdot \xi} \hat{f}(\xi) \Big\}.
\end{align*}
We claim that when $-1 \leq \delta \leq 1$, 
\begin{equation} \label{tuestimate}
\norm{Tu}_{L^2_{\delta}} \leq C \abs{\zeta}^2 \norm{u}_{L^2_{\delta}}.
\end{equation}
For $\delta = 0$ this follows by Fourier multiplier properties since $\abs{\xi_j \xi_k \phi(\xi)} \leq C\abs{\zeta}^2$. For $\delta = 1$ the statement is equivalent with 
\begin{equation*}
\norm{\widehat{Tu}}_{H^1} \leq C \abs{\zeta}^2 \norm{\hat{u}}_{H^1}.
\end{equation*}
One has 
\begin{equation*}
\nabla(\widehat{Tu}(\xi)) = \nabla (\xi_j \xi_k \phi(\xi)) \hat{u}(\xi) + \xi_j \xi_k \phi(\xi) \nabla \hat{u}(\xi).
\end{equation*}
Since $\abs{\nabla (\xi_j \xi_k \phi(\xi))} \leq C \abs{\zeta}$, we have $\norm{\widehat{Tu}}_{H^1} \leq C(\abs{\zeta}^2 + \abs{\zeta}) \norm{\hat{u}}_{H^1} \leq C\abs{\zeta}^2 \norm{\hat{u}}_{H^1}$ using $\abs{\zeta} \geq 1$. This gives \eqref{tuestimate} for $\delta = 1$, duality gives the estimate for $\delta = -1$, and an interpolation gives \eqref{tuestimate} for $-1 \leq \delta \leq 1$.

For $S$ we claim that 
\begin{equation} \label{sfestimate}
\norm{Sf}_{L^2_{\delta}} \leq C \norm{f}_{L^2_{\delta}}
\end{equation}
for $-1 \leq \delta \leq 1$. We may write $Sf = \mF^{-1}\{\psi \hat{f}\}$ where $\psi(\xi) = \frac{\xi_j \xi_k (1-\phi(\xi))}{\abs{\xi}^2 + 2\zeta \cdot \xi}$. If $\abs{\xi} \geq 4\abs{\zeta}$ then $\abs{\abs{\xi}^2 + 2\zeta \cdot \xi} \geq \abs{\xi}^2 - 2\abs{\zeta}\,\abs{\xi} \geq \abs{\xi}^2/2$, so $\abs{\psi} \leq C$ and $\abs{\nabla \psi} \leq C/\abs{\zeta}$. Similar computations as for $T$ imply \eqref{sfestimate} for $-1 \leq \delta \leq 1$. Then \eqref{htwodecomposition}, the $s=0$ estimate for $\Delta_{\zeta}^{-1}$, \eqref{tuestimate}, \eqref{sfestimate}, and the embedding $L^2_{\delta+1} \to L^2_{\delta}$ imply 
\begin{equation*}
\norm{D_j D_k \Delta_{\zeta}^{-1} f}_{L^2_{\delta}} \leq C \abs{\zeta} \norm{f}_{L^2_{\delta+1}}
\end{equation*}
for $-1 < \delta < 0$ and $\abs{\zeta} \geq 1$. This gives the $s=2$ estimate for $\Delta_{\zeta}^{-1}$, and an interpolation gives the estimate for $0 \leq s \leq 2$.
\end{proof}

An easy perturbation argument using Proposition \ref{deltarhoinvestimates} proves Theorem \ref{normestimates} in the case where $W = 0$. The perturbation argument fails in the case where $W$ is nonzero and large. Following Nakamura and Uhlmann \cite{nakamurauhlmann} we will use pseudodifferential operators to conjugate the first order term into a zero order term, so that the perturbation argument can be applied.

First we will write the equation in semiclassical notation. If $\zeta \in \mC^n$ satisfies $\zeta \cdot \zeta = 0$, then we have $\zeta = \mu/h$ with $\mu = \gamma_1 + i\gamma_2$ where $\gamma_1,\gamma_2 \in \mR^n$ satisfy $\abs{\gamma_j} = 1$ and $\gamma_1 \cdot \gamma_2 = 0$, and $h = \sqrt{2}/\abs{\zeta}$ is a small parameter. Note that Proposition \ref{deltarhoinvestimates} gives $\norm{\Delta_{\zeta}^{-1} f}_{H^1_{\delta,h}} \leq Ch \norm{f}_{L^2_{\delta+1}}$, and one also has $\norm{D_{\zeta} f}_{L^2_{\delta}} \leq Ch^{-1} \norm{f}_{H^1_{\delta,h}}$. With this notation,
\begin{align*}
\Delta_{\zeta} &= h^{-2} q(hD), \\
W \cdot D_{\zeta} &= h^{-1} r(x,hD)
\end{align*}
with $q(\xi) = \xi^2 + 2 \mu \cdot \xi$ and $r(x,\xi) = W(x) \cdot (\xi + \mu)$. To take care of the nonsmooth symbol $r$ we will make a $h$-dependent decomposition $r = r^{\sharp} + r^{\flat}$, where $r^{\sharp}$ is a smooth approximation and $r^{\flat}$ is a remainder which will have small norm on suitable spaces when $h$ is small.

Next, we will show that one may conjugate the first order term $r^{\sharp}$ into a zero order term. In the proof we need some facts about the symbol $q$. Since $q(\xi) = (\xi+\gamma_1)^2 - 1 + 2i \gamma_2 \cdot \xi$, this symbol has zero set 
\begin{equation*}
q^{-1}(0) = \{ \xi \in \mR^n \,;\, \abs{\xi+\gamma_1} = 1, \xi \cdot \gamma_2 = 0 \}.
\end{equation*}
For $\varepsilon > 0$ we will consider the neighborhood 
\begin{equation*}
U(\varepsilon) = \{ \xi \in \mR^n \,;\, 1-\varepsilon < \abs{\xi+\gamma_1} < 1+\varepsilon, \abs{\xi \cdot \gamma_2} < \varepsilon \}.
\end{equation*}

\begin{lemma} \label{conjugatingoperators}
Let $0 < \sigma_0 < \sigma < 1/2$ and let $r^{\sharp} \in S^1_{\sigma_0}(\mR^n)$ have the special form $r^{\sharp}(x,\xi) = W^{\sharp}(x) \cdot (\xi + \mu)$, where $W^{\sharp} \in C^{\infty}(\mR^n ; \mC^n)$ satisfies 
\begin{eqnarray*}
 & \abs{\partial^{\alpha} W^{\sharp}(x)} \leq C_{\alpha} h^{-\sigma_0 \abs{\alpha}}, & \\
 & W^{\sharp}(x) = 0 \ \text{for } \abs{x} \geq M. & 
\end{eqnarray*}
Then there exist $a, b, r_0 \in S^0_{\sigma}(\mR^n)$ and $\varepsilon = \varepsilon(\sigma_0,\sigma) > 0$ such that 
\begin{equation} \label{intertwiningeq}
(Q + 2hR^{\sharp})A = BQ + h^{1+\varepsilon} R_0.
\end{equation}
Further, $\br{x} r_0 \in S^0_{\sigma}$ so that $R_0$ is bounded from $L^2_{\delta}$ to $L^2_{\delta+1}$, and for any $s_0, \delta_0 > 0$ there is $h_0 \leq 1$ such that whenever $h \leq h_0$, $A$ and $B$ are bounded and have bounded inverses on $H^s_{\delta,h}$ for $\abs{s} \leq s_0, \abs{\delta} \leq \delta_0$, all with norms bounded uniformly in $h$.
\end{lemma}
\begin{proof}
If $a \in S^0_{\sigma}$ then a direct computation using the special forms of $Q$ and $R^{\sharp}$ implies that 
\begin{equation} \label{firstintertwiningeq}
(Q + 2hR^{\sharp})A = AQ + h \mOp_h(\frac{1}{i} H_q a + 2r^{\sharp}a) + h^2 \mOp_h(\Delta_x a + 2W^{\sharp} \cdot D_x a).
\end{equation}
The last term is in $h^{2-2\sigma} \mOp_h\,S^0_{\sigma}$, and looking at \eqref{intertwiningeq} we would like the middle term to vanish. If $a = e^{i\phi}$ with $\phi \in S^0_{\sigma}$, this would mean that 
\begin{equation*}
(\xi + \mu) \cdot \nabla_x \phi = -r^{\sharp}
\end{equation*}
since $H_q = 2(\xi+\mu) \cdot \nabla_x$. The operator $(\xi + \gamma_1 + i\gamma_2) \cdot \nabla_x$ looks like $\partial_{x_1} + i \partial_{x_2}$ in different coordinates provided that $\xi \in q^{-1}(0)$, but degenerates away from $q^{-1}(0)$. Therefore we will only work in a neighborhood of $q^{-1}(0)$ and introduce a cutoff $\psi(\xi) \in C^{\infty}_c(\mR^n)$ with $\psi = 1$ in $U(\varepsilon/4)$ and $\psi = 0$ outside of $U(\varepsilon/2)$, with $\varepsilon$ as in Lemma \ref{dbarmodified}. This will give a symbol $w(x,\xi)$ with the following properties.

\begin{lemma}
The function 
\begin{equation}
w(x,\xi) = -\frac{1}{2\pi} \int_{\mR^2} \frac{1}{y_1+iy_2} \psi(\xi) r^{\sharp}(x-y_1 (\xi+\gamma_1)-y_2 \gamma_2,\xi) \,dy_1 \,dy_2
\end{equation}
is $C^{\infty}$, solves the equation 
\begin{equation*}
(\xi + \mu) \cdot \nabla_x w = -\psi(\xi)r^{\sharp}(x,\xi),
\end{equation*}
and satisfies the estimates 
\begin{equation} \label{wbounds}
\abs{\partial_x^{\alpha} \partial_{\xi}^{\beta} w(x,\xi)} \leq C_{\alpha \beta} h^{-\sigma_0 \abs{\alpha+\beta}} \br{x}^{\abs{\beta}-1}.
\end{equation}
\end{lemma}
\begin{proof}
This follows from Lemma \ref{dbarmodified} with $V = U(\varepsilon)$, $\gamma_1(\xi) = \xi + \gamma_1$, and $\gamma_2(\xi) = \gamma_2$.
\end{proof}

Note that $w$ is compactly supported in $\xi$ but does not have good behaviour in $x$. To take care of this we will need another cutoff $\chi(x) \in C^{\infty}_c(\mR^n)$ with $\chi = 1$ on $\closure{B(0,M)}$, and we will define 
\begin{equation*}
\phi(x,\xi) = \chi(h^{\theta} x) w(x,\xi)
\end{equation*}
where $\theta = \sigma - \sigma_0$. Then \eqref{wbounds} shows that $\phi$ satisfies 
\begin{equation*}
\abs{\partial_x^{\alpha} \partial_{\xi}^{\beta} \phi(x,\xi)} \leq C_{\alpha \beta} \br{x}^{-1} h^{-\sigma\abs{\alpha+\beta}}.
\end{equation*}
Thus $\phi$ and $\br{x}\phi$ are in $S^0_{\sigma}$, and also $a = e^{i\phi}$ is in $S^0_{\sigma}$. We have 
\begin{multline*}
\frac{1}{i} H_q a + 2r^{\sharp}a = 2(\xi+\mu) \cdot (\nabla_x \phi)e^{i\phi} + 2r^{\sharp} e^{i\phi} \\
 = 2 e^{i\phi} \Big[ (\xi+\mu) \cdot \chi(h^{\theta} x) \nabla_x w + r^{\sharp} + h^{\theta} (\xi+\mu) \cdot w \nabla \chi(h^{\theta} x) \Big].
\end{multline*}
Since $h \leq 1$ we have $r^{\sharp} = \chi(h^{\theta}x) r^{\sharp}$ and 
\begin{equation*}
\frac{1}{i} H_q a + 2r^{\sharp}a = 2 e^{i\phi} \Big[ (1-\psi(\xi)) r^{\sharp} + h^{\theta} (\xi+\mu) \cdot w \nabla \chi(h^{\theta}x) \Big].
\end{equation*}
The second term is compactly supported in $x$ and $\xi$, but the first term is of first order and we are not yet in the situation of \eqref{intertwiningeq}. Here we are saved by the fact that the main operator $Q$ in \eqref{intertwiningeq} is elliptic on the support of $1-\psi(\xi)$, and we may define $b \in S^0_{\sigma}$ by 
\begin{equation*}
b = a + 2h \frac{1-\psi(\xi)}{q(\xi)} e^{i\phi} r^{\sharp}.
\end{equation*}
Then \eqref{firstintertwiningeq} becomes 
\begin{multline*}
(Q + 2hR^{\sharp})A = BQ + h^{1+\theta} \mOp_h(2e^{i\phi} [(\xi+\mu) \cdot w  \nabla \chi(h^{\theta}x)]) \\
 + h^{2-2\sigma} \mOp_h(h^{2\sigma} [\Delta_x a + 2W^{\sharp} \cdot D_x a]).
\end{multline*}
Choosing $\varepsilon = \min \{ \theta, 1-2\sigma \}$ gives \eqref{intertwiningeq} with $r_0 \in S^0_{\sigma}$, and one even has $\br{x} r_0 \in S^0_{\sigma}$. Proposition \ref{weightedpsdoestimates} gives that $R_0$ is bounded from $L^2_{\delta}$ to $L^2_{\delta+1}$ with norm bounded uniformly in $h$.

It remains to show that $A$ and $B$ are bounded and invertible on $H^s_{\delta,h}$ for $\abs{s} \leq s_0, \abs{\delta} \leq \delta_0$ when $h$ is small enough, and that all norms are bounded uniformly in $h$. The boundedness of $A$ and $B$ follows from Proposition \ref{weightedpsdoestimates}. To show invertibility we note that $1/a = e^{-i\phi}$ is in $S^0_{\sigma}$ and 
\begin{equation*}
\mOp_h(a) \mOp_h(1/a) = I + h^{1-2\sigma} \mOp_h(m)
\end{equation*}
where $m \in S^0_{\sigma}$ by Proposition \ref{psdoproperties}. Since $\mOp_h(m)$ has bounded norm on $H^s_{\delta,h}$ for $\abs{s} \leq s_0$, $\abs{\delta} \leq \delta_0$, the operator $I + h^{1-2\sigma} \mOp_h(m)$ is invertible on these spaces if $h$ is small enough. Then also $A$ is invertible with norm of the inverse uniformly bounded in $h$. The same applies to $B$ since $b = a + h \mOp_h\,S^{-1}_{\sigma}$.
\end{proof}

\begin{remark}
Lemma \ref{conjugatingoperators} is a global nonsmooth version of the pseudodifferential conjugation technique in \cite{nakamurauhlmann} (see also \cite{nakamurauhlmannerratum}). Similar ideas have been used in inverse scattering \cite{eskinralston}, \cite{isozakinakazawauhlmann}, nonlinear Schr\"odinger equations \cite{takeuchi}, \cite{kenigponcevega} and periodic Schr\"odinger operators \cite{sobolev}.

The problem in extending the method to the global case is seen in \eqref{wbounds}, where the derivatives in $\xi$ of the symbol grow in $x$. This behaviour leads to poor global properties. A solution, presented in \cite{takeuchi} and \cite{kenigponcevega}, was to multiply a symbol $a$ satisfying $\abs{\partial_x^{\alpha} \partial_{\xi}^{\beta} a(x,\xi)} \leq C_{\alpha \beta} \br{x}^{\abs{\beta}} \br{\xi}^{-\abs{\beta}}$ by a cutoff $\chi(R_0 \br{x}/\br{\xi})$. The new symbol is of type $(0,0)$, hence bounded on $L^2$, and the error term which appears in the equation because of this modification is of lower order. One can even invert related operators on $L^2$ by adjusting the parameter $R_0$.

In the present case there is an additional parameter $h$ which may be taken arbitrarily small, and additionally $\abs{\xi} \leq C h^{-1}$ in the support of $w(x,h\xi)$. Therefore our cutoff has the simpler form $\chi(h^{\theta} x)$, where $\theta$ is chosen so that the new symbol falls into $S^0_{\sigma}$ with $\sigma < 1/2$.
\end{remark}

We proceed to prove the main norm estimates.

\begin{proof}
(of Theorem \ref{normestimates}) The proof is given in three steps.

\medskip

\noindent \emph{Step 1}: A decomposition

\medskip

Let $\varphi \in C_c^{\infty}(\mR^n)$, $\varphi \geq 0$, $\varphi = 1$ for $\abs{x} \leq 1/2$, and $\varphi = 0$ for $\abs{x} \geq 1$. Let also $\varphi_{\varepsilon}(x) = \varepsilon^{-n} \varphi(x/\varepsilon)$ be the usual mollifier. We will use the decomposition (as in \cite{salo}) 
\begin{equation} \label{wdecomposition}
W = W^{\sharp} + W^{\flat}
\end{equation}
where $W^{\sharp} = W \ast \varphi_{\varepsilon}$ is a smooth approximation of $W$, and we make the specific choice 
\begin{equation*}
\varepsilon = h^{\sigma_0}
\end{equation*}
where $0 < \sigma_0 < 1/2$. Then $W^{\flat}$ is a small remainder term, and one has 
\begin{eqnarray*}
 & \abs{\partial^{\alpha} W^{\sharp}(x)} \leq C_{\alpha} h^{-\sigma_0 \abs{\alpha}}, & \\
 & \norm{W^{\flat}}_{L^{\infty}} \to 0 \text{ as } h \to 0, &
\end{eqnarray*}
the second estimate by the continuity of $W$.

\medskip

\noindent \emph{Step 2}: Existence

\medskip

Using the decomposition \eqref{wdecomposition}, we write the equation \eqref{mainequation} as 
\begin{equation} \label{maineqdecomposition}
(\Delta_{\zeta} + 2W^{\sharp} \cdot D_{\zeta} + 2W^{\flat} \cdot D_{\zeta} + q) u = f.
\end{equation}
Choosing $\sigma$ with $\sigma_0 < \sigma < 1/2$ and $s_0 = \delta_0 = 2$, Lemma \ref{conjugatingoperators} gives $h_0 \leq 1$ and $a,b,r_0 \in S^0_{\sigma}$ with 
\begin{equation} \label{maineqconjugation}
(\Delta_{\zeta} + 2W^{\sharp} \cdot D_{\zeta}) A = B \Delta_{\zeta} + h^{-1+\varepsilon} R_0.
\end{equation}
We will assume $h \leq h_0$, so $A$ and $B$ will be invertible. We look for a solution of \eqref{maineqdecomposition} of the form $u = \Delta_{\zeta}^{-1} v$ for $v \in L^2_{\delta+1}$. Then $u = A A^{-1} \Delta_{\zeta}^{-1} v$, and inserting this in \eqref{maineqdecomposition} and using \eqref{maineqconjugation} gives 
\begin{equation*}
(B \Delta_{\zeta} A^{-1} \Delta_{\zeta}^{-1} + h^{-1+\varepsilon} R_0 A^{-1} \Delta_{\zeta}^{-1} + 2W^{\flat} \cdot D_{\zeta} \Delta_{\zeta}^{-1} + q \Delta_{\zeta}^{-1}) v = f.
\end{equation*}
We will show that in the operator on the left, the last three terms are small perturbations of the first term when $h$ is small.

Consider the operator 
\begin{equation*}
M = B \Delta_{\zeta} A^{-1} \Delta_{\zeta}^{-1}.
\end{equation*}
From \eqref{maineqconjugation} we get 
\begin{equation*}
M = I + 2 W^{\sharp} \cdot D_{\zeta} \Delta_{\zeta}^{-1} - h^{-1+\varepsilon} R_0 A^{-1} \Delta_{\zeta}^{-1},
\end{equation*}
and therefore $M$ is bounded on $L^2_{\delta+1}$ with norm bounded uniformly in $h$. It is easy to see that $M$ has the inverse 
\begin{equation*}
M^{-1} = \Delta_{\zeta} A \Delta_{\zeta}^{-1} B^{-1}.
\end{equation*}
Similarly from \eqref{maineqconjugation} we obtain 
\begin{equation*}
M^{-1} = I - 2 W^{\sharp} \cdot D_{\zeta} A \Delta_{\zeta}^{-1} B^{-1} + h^{-1+\varepsilon} R_0 \Delta_{\zeta}^{-1} B^{-1}
\end{equation*}
which is again bounded on $L^2_{\delta+1}$, with norm bounded uniformly in $h$. Also, using the mapping properties of the related operators and the decay of $\norm{W^{\flat}}_{L^{\infty}}$, we have 
\begin{equation*}
\norm{h^{-1+\varepsilon} R_0 A^{-1} \Delta_{\zeta}^{-1} + 2W^{\flat} \cdot D_{\zeta} \Delta_{\zeta}^{-1} + q \Delta_{\zeta}^{-1}}_{L^2_{\delta+1} \to L^2_{\delta+1}} = o(1)
\end{equation*}
as $h \to 0$. Then we obtain a solution $u$ of \eqref{mainequation} in the form 
\begin{equation*}
u = \Delta_{\zeta}^{-1} M^{-1} (I + h^{-1+\varepsilon} R_0 A^{-1} \Delta_{\zeta}^{-1} M^{-1} + 2W^{\flat} \cdot D_{\zeta} \Delta_{\zeta}^{-1} M^{-1} + q \Delta_{\zeta}^{-1} M^{-1})^{-1} f.
\end{equation*}
Thus $u = \Delta_{\zeta}^{-1} v$ with $\norm{v}_{L^2_{\delta+1}} \leq C \norm{f}_{L^2_{\delta+1}}$. The norm estimates for $u$ follow from Proposition \ref{deltarhoinvestimates}.

\medskip

\noindent \emph{Step 3}: Uniqueness

\medskip

It is enough to show that if $u \in H^1_{\delta}$ satisfies 
\begin{equation} \label{uniquenesseq}
(\Delta_{\zeta} + 2W^{\sharp} \cdot D_{\zeta} + 2W^{\flat} \cdot D_{\zeta} + q)u = 0,
\end{equation}
then $u = 0$. We use Lemma \ref{conjugatingoperators}. It follows that $u = A v$ for $v = A^{-1} u \in H^1_{\delta}$, so that $v$ satisfies 
\begin{equation*}
(B \Delta_{\zeta} + h^{-1+\varepsilon} R_0 + 2W^{\flat} \cdot D_{\zeta} A + q A) v = 0.
\end{equation*}
Applying $B^{-1}$ from the left we get 
\begin{equation} \label{deltarhoveq}
\Delta_{\zeta} v = -(h^{-1+\varepsilon} B^{-1} R_0 + 2 B^{-1} W^{\flat} \cdot D_{\zeta} A + B^{-1} q A)v.
\end{equation}

The right hand side of \eqref{deltarhoveq} is in $L^2_{\delta+1}$ since $R_0$ and $q A$ map $L^2_{\delta}$ to $L^2_{\delta+1}$ and $W^{\flat} \cdot D_{\zeta} A$ maps $H^1_{\delta}$ to $L^2_{\delta+1}$. We are now in the situation of Proposition \ref{deltarhoinvestimates}, and using the $H^1_{\delta}$ estimate of that Proposition to \eqref{deltarhoveq} implies 
\begin{equation*}
\norm{v}_{H^1_{\delta}} \leq C( \norm{ (h^{-1+\varepsilon} B^{-1} R_0 + B^{-1} q A) v }_{L^2_{\delta+1}} + \norm{2 B^{-1} W^{\flat} \cdot D_{\zeta} A v}_{L^2_{\delta+1}} ).
\end{equation*}
The first term on the right is bounded by $C \abs{\zeta} \norm{v}_{L^2_{\delta}}$ and the second term is bounded by $C \norm{\br{x} W^{\flat}}_{L^{\infty}} (\norm{v}_{H^1_{\delta}} + \abs{\zeta} \norm{v}_{L^2_{\delta}})$, with $C$ independent of $\zeta$. Choosing $\abs{\zeta}$ so large that $\norm{\br{x} W^{\flat}}_{L^{\infty}} \leq 1/(2C)$, we have that the coefficient of $\norm{v}_{H^1_{\delta}}$ on the right is $\leq 1/2$, and we may move this term to the left. We are left with the estimate 
\begin{equation} \label{vh1estimate}
\norm{v}_{H^1_{\delta}} \leq C \abs{\zeta} \norm{v}_{L^2_{\delta}}
\end{equation}
with $C$ independent of $\zeta$.

Finally, we use the $L^2_{\delta}$ estimate of Proposition \ref{deltarhoinvestimates} to \eqref{deltarhoveq}. This implies 
\begin{equation*}
\norm{v}_{L^2_{\delta}} \leq \frac{C}{\abs{\zeta}} ( \norm{ (h^{-1+\varepsilon} B^{-1} R_0 + B^{-1} q A) v }_{L^2_{\delta+1}} + \norm{ 2 B^{-1} W^{\flat} \cdot D_{\zeta} A v }_{L^2_{\delta+1}} )
\end{equation*}
The first term inside the parentheses is $\leq C \abs{\zeta}^{1-\varepsilon} \norm{v}_{L^2_{\delta}}$, and the second is $\leq C \norm{\br{x} W^{\flat}}_{L^{\infty}} (\norm{v}_{H^1_{\delta}} + \abs{\zeta}\,\norm{v}_{L^2_{\delta}})$. Using \eqref{vh1estimate} this becomes 
\begin{equation*}
\norm{v}_{L^2_{\delta}} \leq C (\abs{\zeta}^{-\varepsilon} + \norm{\br{x} W^{\flat}}_{L^{\infty}}) \norm{v}_{L^2_{\delta}}
\end{equation*}
where $C$ is independent of $\zeta$. Choosing $\abs{\zeta}$ large enough we obtain $\norm{v}_{L^2_{\delta}} \leq \frac{1}{2} \norm{v}_{L^2_{\delta}}$, implying $v = 0$ and also $u = 0$.
\end{proof}

\section{Equivalent problems}

In this section we show that complex geometrical optics solutions for the magnetic Schr\"odinger equation can be characterized in several different ways. The treatment is almost completely analogous with \cite{nachman}. We begin by stating the main result and explain the notation later as we go along.

\begin{prop} \label{equivalentproblems}
Let $\Omega \subseteq \mR^n$ be a bounded domain with $C^{1,1}$ boundary. Suppose $W \in L^{\infty}_{\Omega}(\mR^n ; \mC^n)$ with $D \cdot W \in L^{\infty}$, and suppose $q \in L^{\infty}_{\Omega}(\mR^n ; \mC)$. Also suppose that $0$ is not a Dirichlet eigenvalue of $H_{W,q}$ in $\Omega$. Let $\zeta \in \mC^n$ with $\zeta^2 = 0$, and consider the following four problems:
\begin{align*}
\text{(DE)} & \left\{ \begin{array}{l}
H_{W,q} u = 0 \text{ in } \mR^n \\[2pt] 
u = e^{i\zeta \cdot x}(1+\omega) \text{ with } \omega \in \Delta_{\zeta}^{-1} L^2_{\Omega}, 
\end{array} \right. \\ %\\[12pt]
\text{(IE)\,} & \left\{ \begin{array}{l} 
u + G_{\zeta} \ast (2W \cdot Du + (W^2 + D \cdot W + q)u) = e^{i\zeta \cdot x} \text{ in } \mR^n \\[2pt] 
u \in H^1_{\mathrm{loc}}(\mR^n), 
\end{array} \right. \\ %\\[12pt]
\text{(EP)} & \left\{ \begin{array}{rl} 
\text{i)} & \Delta u = 0 \text{ in } \Omega' \\[2pt]
\text{ii)} & u \in H^2(\Omega_R') \text{ for any } R > R_0 \\[2pt]
\text{iii)} & u \text{ satisfies \eqref{radiationcondition} for almost every } x \in \mR^n \\[2pt]
\text{iv)} & \frac{\partial u}{\partial \nu_+} = \Lambda_{W,q}(u_+) \text{ on } \partial \Omega,
\end{array} \right. \\ %\\[32pt]
\text{(BE)} & \left\{ \begin{array}{l} 
(\frac{1}{2}I + S_{\zeta} \Lambda_{W,q} - B_{\zeta}) f = e^{i\zeta \cdot x} \text{ on } \partial \Omega \\[2pt] 
f \in H^{3/2}(\partial \Omega). 
\end{array} \right.
\end{align*}
Then all these problems are equivalent, in the sense that if a solution exists (is unique) for one problem, then a solution exists (is unique) for all the problems. If $u$ is a solution of (DE), then $u$ solves (IE), $u|_{\Omega'}$ solves (EP), and $u|_{\partial \Omega}$ solves (BE).
\end{prop}

\begin{remark}
The assumption that $0$ is not an eigenvalue is for simplicity. With appropriate changes, similar results are valid also when $0$ is an eigenvalue.
\end{remark}

The first step is to show the equivalence of the differential equation (DE) and the integral equation (IE). This involves the Green function $G_{\zeta}$, defined by 
\begin{equation*}
G_{\zeta} = e^{i\zeta \cdot x} g_{\zeta}
\end{equation*}
where $g_{\zeta}$ is the tempered distribution such that $\Delta_{\zeta}^{-1} f = g_{\zeta} \ast f$ for $f$ in the Schwartz class ($g_{\zeta}$ exists since $\Delta_{\zeta}^{-1}$ is translation invariant). Then 
\begin{equation*}
\Delta G_{\zeta} = \zeta^2 G_{\zeta} + 2 e^{i\zeta \cdot x} \zeta \cdot D g_{\zeta} + e^{i\zeta \cdot x} \Delta g_{\zeta} = e^{i\zeta \cdot x} \Delta_{\zeta} g_{\zeta} = \delta_0
\end{equation*}
where $\delta_0$ is the Dirac measure at $0$. Consequently $G_{\zeta} = G_0 + H_{\zeta}$ where $G_0(x) = c_n \abs{x}^{2-n}$ is the usual fundamental solution of $\Delta$, and $H_{\zeta}$ is a global harmonic function (one has $c_n = \frac{1}{n(n-2)\alpha(n)}$ where $\alpha(n)$ is the volume of the $n$-dimensional unit ball).

We note that the left hand side of (IE) is well defined for any $u \in H^1_{\mathrm{loc}}(\mR^n)$, since then $2W \cdot Du + (W^2 + D \cdot W + q)u \in L^2_{\Omega}(\mR^n)$ is a compactly supported distribution. Also note that $G_{\zeta} \ast f = e^{i\zeta \cdot x} \Delta_{\zeta}^{-1} e^{-i\zeta \cdot x} f$ whenever $f \in L^2_c(\mR^n)$.

\begin{lemma} \label{deieequivalence}
Assume the conditions of Proposition \ref{equivalentproblems}. Then, $u$ is a solution of (DE) if and only if $u$ is a solution of (IE). Also, a solution of (DE) is unique if and only if a solution of (IE) is unique.
\end{lemma}
\begin{proof}
Let first $u = e^{i\zeta \cdot x}(1+\omega)$ solve (DE) where $\omega = \Delta_{\zeta}^{-1} f$ with $f \in L^2_{\Omega}$. Clearly $u \in H^1_{\mathrm{loc}}(\mR^n)$, and $H_{W,q} u = 0$ implies 
\begin{equation*}
(\Delta_{\zeta} + 2 W \cdot D_{\zeta} + (W^2 + D \cdot W + q))(1 + \Delta_{\zeta}^{-1} f) = 0.
\end{equation*}
We have $\Delta_{\zeta}(1 + \Delta_{\zeta}^{-1} f) = f$. Now applying $\Delta_{\zeta}^{-1}$ to both sides, which is allowed since the left hand side is in $L^2_{\Omega}$, gives 
\begin{equation*}
\omega + \Delta_{\zeta}^{-1} (2W \cdot D_{\zeta}(1 + \omega) + (W^2 + D \cdot W + q)(1 + \omega)) = 0.
\end{equation*}
We obtain (IE) by adding the constant one to both sides and multiplying by $e^{i\zeta \cdot x}$.

For the converse, suppose $u$ solves (IE), and write $u = e^{i\zeta \cdot x} u_0$. Then $u_0$ solves $u_0 + \Delta_{\zeta}^{-1}(2W \cdot D_{\zeta} u_0 + (W^2 + D \cdot W + q) u_0) = 1$. Applying $\Delta_{\zeta}$ to both sides gives $H_{W,q} u = 0$. Also, one sees that $u_0 - 1 = \Delta_{\zeta}^{-1} f$ for $f \in L^2_{\Omega}(\mR^n)$.

The uniqueness part is obtained just by noting that if $u_1$ and $u_2$ solve (DE) then $u_1$ and $u_2$ solve (IE), and vice versa.
\end{proof}

Next we show that (IE) and the exterior problem (EP) are equivalent. We have used the notation $\Omega' = \mR^n \smallsetminus \closure{\Omega}$ and $\Omega'_R = B(0,R) \smallsetminus \closure{\Omega}$, where $R > R_0$ and $\closure{\Omega} \subseteq B(0,R_0)$. We write $u_{+}$ (resp. $u_{-}$) for the restriction of $u$ to $\partial \Omega$ from the exterior (resp. interior), and $\frac{\partial u}{\partial \nu_{+}}$ (resp. $\frac{\partial u}{\partial \nu_{-}}$) for the value of $\nabla u \cdot \nu$ on $\partial \Omega$ from the exterior (resp. interior), where $\nu$ is the outer unit normal to $\partial \Omega$. We also write $G_{\zeta}(x,y) = G_{\zeta}(x-y)$.

A main point will be that a solution $u$ of (IE) satisfies the radiation condition 
\begin{equation} \label{radiationcondition}
\int_{\abs{y} = R} \Big( G_{\zeta}(x,y)\frac{\partial u}{\partial \nu}(y) - u(y) \frac{\partial G_{\zeta}(x,y)}{\partial \nu(y)} \Big) \,dS(y) \to e^{i\zeta \cdot x}
\end{equation}
for a.e.~$x \in \mR^n$ as $R \to \infty$. For applications of Green's identity below, we define a smooth approximation of $G_{\zeta}$ by $G_{\zeta}^{\varepsilon} = G_0^{\varepsilon} + H_{\zeta}$, where 
\begin{equation*}
G_0^{\varepsilon}(x) = c_n (\varepsilon^2 + \abs{x}^2)^{\frac{2-n}{2}}.
\end{equation*}
Note that $\Delta G_{\zeta}^{\varepsilon}(x) = \varepsilon^{-n}\varphi(x/\varepsilon)$ where 
\begin{equation*}
\varphi(x) = \frac{1}{\alpha(n)} (1+\abs{x}^2)^{-\frac{n+2}{2}}
\end{equation*}
and $\int \varphi(x) \,dx = 1$. Thus $\Delta G_{\zeta}^{\varepsilon}$ is an approximation of the identity.

Before showing the equivalence of (IE) and (EP) we need a lemma on regularity properties of solutions of $H_{W,q} u = 0$ and of $\Lambda_{W,q}$.

\begin{lemma} \label{regulardnmap}
Under the conditions of Proposition \ref{equivalentproblems}, the operator $P_{W,q}$, which maps $f \in H^{3/2}(\Omega)$ to the solution $u$ of $H_{W,q} u = 0$ in $\Omega$ with $u|_{\partial \Omega} = f$, is bounded $H^{3/2}(\partial \Omega) \to H^2(\Omega)$. Further, one has $\Lambda_{W,q}: H^{3/2}(\partial \Omega) \to H^{1/2}(\partial \Omega)$, and 
\begin{equation*}
\Lambda_{W,q} f = \frac{\partial u}{\partial \nu} \Big|_{\partial \Omega}.
\end{equation*}
\end{lemma}
\begin{proof}
The operator $H_{W,q}$, written in nondivergence form, satisfies the assumptions of \cite[Theorem 8.12]{gilbargtrudinger} (the theorem is given for $C^2$ domains but the result holds with the same proof for $C^{1,1}$ domains). This shows that $u$ is in $H^2(\Omega)$ if $f \in H^{3/2}(\partial \Omega)$, and that the solution operator $P_{W,q}$ is bounded.

For the second part, we claim that if $W \in L^n_{\Omega}(\mR^n ; \mC^n)$ and $D \cdot W \in L^{n/2}(\mR^n ; \mC)$, then for any $v \in W^{1,n/(n-1)}(\Omega)$ one has 
\begin{equation} \label{wnormalweakzero}
\int_{\Omega} (W \cdot Dv + (D \cdot W)v) \,dx = 0.
\end{equation}
This statement means that $W \cdot \nu = 0$ on $\partial \Omega$, in a certain weak sense. The expression is well defined since $v \in L^{n/(n-2)}$ by Sobolev embedding. We take $W_j \in C^{\infty}_c(\mR^n ; \mC^n)$ to be convolution approximations of $W$ so that $W_j \to W$ in $L^n$ and $D \cdot W_j \to D \cdot W$ in $L^{n/2}$, and we take an extension of $v$ in $W^{1,n/(n-1)}(\mR^n)$. If the supports of $W_j$ and $W$ are contained in $B(0,R)$, then 
\begin{align*}
\int_{\Omega} (W \cdot Dv + (D \cdot W)v) \,dx &= \lim_{j \to \infty} \int_{B(0,R)} (W_j \cdot Dv + (D \cdot W_j)v) \,dx \\
 &= \lim_{j \to \infty} \frac{1}{i} \int_{\partial B(0,R)} (W_j \cdot \nu) v \,dS = 0.
\end{align*}
Now let $f, g \in H^{3/2}(\Omega)$ and let $u_f = P_{W,q} f$ and $e_g \in H^2(\Omega)$ with $e_g|_{\partial \Omega} = g$. An integration by parts gives 
\begin{equation} \label{unormalweak}
\langle \frac{\partial u_f}{\partial \nu}\Big|_{\partial \Omega}, g \rangle = \int_{\Omega} (\nabla u_f \cdot \nabla e_g + (2W \cdot D u_f + (W^2 + D \cdot W + q) u_f) e_g) \,dx.
\end{equation}
Now $u_f e_g \in W^{2,1}(\Omega) \subseteq W^{1,n/(n-1)}(\Omega)$. Using \eqref{wnormalweakzero} with $v = u_f e_g$ and substituting this to \eqref{unormalweak} gives $\frac{\partial u_f}{\partial \nu}\Big|_{\partial \Omega} = \Lambda_{W,q} f$. This also shows that $\Lambda_{W,q}$ is bounded $H^{3/2}(\partial \Omega) \to H^{1/2}(\partial \Omega)$.
\end{proof}

\begin{lemma} \label{ieepequivalence}
Assume the conditions of Proposition \ref{equivalentproblems}. Then, if $u$ is a solution of (IE), then $u|_{\Omega'}$ is a solution of (EP). Conversely, if $u$ is a solution of (EP), then there is a unique extension $\tilde{u}$ of $u$ to $\mR^n$ so that $\tilde{u}$ is a solution of (IE). Also, a solution of (IE) is unique if and only if a solution of (EP) is unique.
\end{lemma}
\begin{proof}
Suppose $u$ solves (IE). By Lemma \ref{deieequivalence} we have $H_{W,q} u = 0$ and $u \in H^2_{\delta}(\mR^n)$, which shows (EP) i)-ii). To prove iii) fix $x \in \mR^n$ and let $R > \abs{x}$ and $R > R_0$, and write 
\begin{multline*}
\int_{\abs{y} = R} \Big( G_{\zeta}^{\varepsilon}(x,y)\frac{\partial u}{\partial \nu}(y) - u(y) \frac{\partial G_{\zeta}^{\varepsilon}(x,y)}{\partial \nu(y)} \Big) \,dS(y) \\
 = - \int_{B(0,R)} (G_{\zeta}^{\varepsilon}(x,y) \Delta u(y) - u(y) \Delta_y G_{\zeta}^{\varepsilon}(x,y)) \,dy \\
 = \int_{B(0,R)} u \Delta_y G_{\zeta}^{\varepsilon}(x,y) \,dy + \int_{B(0,R)} G_{\zeta}^{\varepsilon}(x,y) (2W \cdot Du + (W^2 + D \cdot W + q)u) \,dy \\
 = (\Delta G_{\zeta}^{\varepsilon} \ast u\chi_{B(0,R)})(x) + (G_{\zeta}^{\varepsilon} \ast (2W \cdot Du + (W^2 + D \cdot W + q)u))(x)
\end{multline*}
since $W$ and $q$ have their supports inside $B(0,R)$.

As $\varepsilon \to 0$, the first term on the right converges to $u(x)$ outside a set of measure zero (this set depends on $R$, but one may take the union of such sets for countably many $R$). The second term on the right converges to $(G_{\zeta} \ast (2W \cdot Du + (W^2 + D \cdot W + q)u))(x)$ for a.e.~$x \in \mR^n$ by dominated convergence, since $G_{\zeta} \in L^1_{\mathrm{loc}}$ and the other function is in $L^2_c$. Since $\abs{x} < R$ the boundary integrals present no problem and one may replace $G_{\zeta}^{\varepsilon}$ by $G_{\zeta}$. We obtain for a.e.~$x$ 
\begin{multline} \label{radiationie}
\lim_{R \to \infty} \int_{\abs{y} = R} \Big( G_{\zeta}(x,y)\frac{\partial u}{\partial \nu}(y) - u(y) \frac{\partial G_{\zeta}(x,y)}{\partial \nu(y)} \Big) \,dS(y) \\
 = u(x) + (G_{\zeta} \ast (2W \cdot Du + (W^2 + D \cdot W + q)u))(x).
\end{multline}
Since $u$ satisfies (IE) we obtain (EP) iii). Further, since $u \in H^2_{\mathrm{loc}}$ and $0$ is not a Dirichlet eigenvalue of $H_{W,q}$ in $\Omega$, Lemma \ref{regulardnmap} gives 
\begin{equation*}
\frac{\partial u}{\partial \nu_{+}} = \frac{\partial u}{\partial \nu_{-}} = \Lambda_{W,q} u_{-} = \Lambda_{W,q} u_{+}
\end{equation*}
which is (EP) iv).

Let now $u$ solve (EP). We use Lemma \ref{regulardnmap} and let $v = P_{W,q} u_+ \in H^2(\Omega)$, and we define $\tilde{u}(x) = u(x)$ for $x \in \Omega'$ and $\tilde{u}(x) = v(x)$ for $x \in \Omega$. Now $\tilde{u}_{-} = v_{-}  = u_{+} = \tilde{u}_{+}$ and 
\begin{equation*}
\frac{\partial \tilde{u}}{\partial \nu_{-}} = \frac{\partial v}{\partial \nu_{-}} = \Lambda_{W,q} v_{-} = \Lambda_{W,q} u_{+} = \frac{\partial u}{\partial \nu_{+}} = \frac{\partial \tilde{u}}{\partial \nu_{+}}
\end{equation*}
by Lemma \ref{regulardnmap} and (EP) iv). This shows that $\tilde{u} \in H^2_{\mathrm{loc}}(\mR^n)$. By (EP) i) we have $H_{W,q} \tilde{u} = 0$ in $\mR^n$, and then the computation above leads to \eqref{radiationie} with $u$ replaced by $\tilde{u}$. The condition (EP) iii) shows that $\tilde{u}$ solves (IE).

The uniqueness part follows from the facts that if $u_1$ and $u_2$ solve (IE) then $u_1|_{\Omega'}$ and $u_2|_{\Omega'}$ solve (EP), and if $u_1$ and $u_2$ solve (EP) then $\tilde{u}_1$ and $\tilde{u}_2$ solve (IE).
\end{proof}

The final equivalence will be between (EP) and the boundary integral equation (BE). Here we need the layer potentials depending on $\zeta$, defined in terms of the Green function $G_{\zeta}$. The single layer potential $S_{\zeta}$, double layer potential $D_{\zeta}$, and boundary layer potential $B_{\zeta}$ are defined by 
\begin{align*}
S_{\zeta} f(x) &= \int_{\partial \Omega} G_{\zeta}(x,y) f(y) \,dS(y) \quad (x \in \mR^n \smallsetminus \partial \Omega), \\
D_{\zeta} f(x) &= \int_{\partial \Omega} \frac{\partial G_{\zeta}(x,y)}{\partial \nu(y)} f(y) \,dS(y) \quad (x \in \mR^n \smallsetminus \partial \Omega), \\
B_{\zeta} f(x) &= \int_{\partial \Omega} \frac{\partial G_{\zeta}(x,y)}{\partial \nu(y)} f(y) \,dS(y) \quad (x \in \partial \Omega).
\end{align*}
Since $\partial \Omega$ is $C^{1,1}$ one does not need a principal value in the definition of $B_{\zeta}$. The operators have the following properties, given in \cite{nachman}.

\begin{prop} \label{layerpotentials}
Let $\Omega \subseteq \mR^n$, $n \geq 3$, be a bounded domain with $C^{1,1}$ boundary, and suppose that $\closure{\Omega} \subseteq B(0,R_0)$.
\begin{enumerate}
\item[(a)] 
Let $f \in H^{1/2}(\partial \Omega)$ and $u = S_{\zeta} f$. Then $u$ is in $H^2(\Omega)$ and $H^2(\Omega_R')$ for any $R > R_0$, and $\Delta u = 0$ in $\mR^n \smallsetminus \partial \Omega$. If $R > R_0$, then $u$ satisfies the radiation condition 
\begin{equation} \label{strongradiationcondition}
\int_{\abs{y} = R} \Big( G_{\zeta}(x,y)\frac{\partial u}{\partial \nu}(y) - u(y) \frac{\partial G_{\zeta}(x,y)}{\partial \nu(y)} \Big) \,dS(y) = 0
\end{equation}
for almost every $x$ with $\abs{x} < R$.
\item[(b)] 
Let $f \in H^{3/2}(\partial \Omega)$ and $v = D_{\zeta} f$. Then $v$ has the properties listed in (a).
\item[(c)] 
In the situation of (a), one has $u_{-} = u_{+}$ on $\partial \Omega$, in the sense of $H^{3/2}(\partial \Omega)$ as well as nontangential convergence a.e.~on $\partial \Omega$. We will write $u = S_{\zeta} f$ on $\partial \Omega$. The map $f \mapsto S_{\zeta} f$ is bounded $H^s(\partial \Omega) \to H^{s+1}(\partial \Omega)$ for $0 \leq s \leq 1$, and one has 
\begin{equation} \label{singlelayerjump}
\frac{\partial u}{\partial \nu_{-}} - \frac{\partial u}{\partial \nu_{+}} = f \quad \text{on } \partial \Omega.
\end{equation}
\item[(d)] 
In the situation of (b), one has 
\begin{equation} \label{doublelayerjump}
v_{\pm} = \pm \frac{1}{2} f + B_{\zeta} f \quad \text{on } \partial \Omega,
\end{equation}
in the sense of $H^{3/2}(\partial \Omega)$ as well as nontangential convergence.
\item[(e)] 
The map $B_{\zeta}$ is bounded $H^s(\partial \Omega) \to H^s(\partial \Omega)$ for $0 \leq s \leq 3/2$.
\end{enumerate}
\end{prop}

\begin{lemma} \label{epbeequivalence}
Assume the conditions of Proposition \ref{equivalentproblems}. Then, if $u$ is a solution of (EP), then $f = u|_{\partial \Omega}$ is a solution of (BE). Conversely, if $f$ is a solution of (BE), then 
\begin{equation} \label{epsolutiondef}
u = e^{i\zeta \cdot x} - S_{\zeta} \Lambda_{W,q} f + D_{\zeta} f
\end{equation}
is a solution of (EP), with $u_{+} = f$. Also, a solution of (EP) is unique if and only if a solution of (BE) is unique.
\end{lemma}
\begin{proof}
Suppose $u$ solves (EP). We let $f = u_{+}$ on $\partial \Omega$. Then $f \in H^{3/2}(\partial \Omega)$. If $x \in \Omega'$ and $R > \abs{x}$, we have 
\begin{multline*}
- \int_{\Omega_R'} (G_{\zeta}^{\varepsilon}(x,y) \Delta u(y) - u(y) \Delta_y G_{\zeta}^{\varepsilon}(x,y)) \,dy \\
 = \Big( \int_{\abs{y} = R} - \int_{\partial \Omega} \Big) \Big( G_{\zeta}^{\varepsilon}(x,y)\frac{\partial u}{\partial \nu}(y) - u(y) \frac{\partial G_{\zeta}^{\varepsilon}(x,y)}{\partial \nu(y)} \Big) \,dS(y).
\end{multline*}
Letting $\varepsilon \to 0$ and using (EP) i) we obtain 
\begin{equation} \label{intermediateformula}
u(x) = \int_{\abs{y} = R} \Big( G_{\zeta}(x,y)\frac{\partial u}{\partial \nu}(y) - u(y) \frac{\partial G_{\zeta}(x,y)}{\partial \nu(y)} \Big) \,dS(y) - S_{\zeta} \Big( \frac{\partial u}{\partial \nu_{+}} \Big)(x) + D_{\zeta} (u_{+})(x)
\end{equation}
for a.e.~$x$ in $\Omega'$. We let $R \to \infty$, use (EP) iii)-iv), and then let $x \to \partial \Omega$ nontangentially and use Proposition \ref{layerpotentials} (d), which gives that $f = u_+$ satisfies (BE).

Conversely, suppose $f$ satisfies (BE) and define $u$ by \eqref{epsolutiondef} in $\Omega'$. Then $u$ satisfies (EP) i)-iii) by Proposition \ref{layerpotentials} (it is an easy calculation that the left hand side of \eqref{strongradiationcondition} equals $e^{i\zeta \cdot x}$ if $u = e^{i\zeta \cdot x}$). We need to show (EP) iv). First note that by Proposition \ref{layerpotentials} (d),
\begin{equation*}
u_{+} = e^{i\zeta \cdot x} - S_{\zeta} \Lambda_{W,q} f + \frac{1}{2} f + B_{\zeta} f \quad \text{on } \partial \Omega,
\end{equation*}
which gives $u_{+} = f$ using (BE). The formula \eqref{intermediateformula} holds for $u$ with the same proof, and as $R \to \infty$ we obtain 
\begin{equation*}
u(x) = e^{i\zeta \cdot x} - S_{\zeta} \Big( \frac{\partial u}{\partial \nu_{+}} \Big)(x) + D_{\zeta} f(x)
\end{equation*}
a.e.~in $\Omega'$. Comparing with \eqref{epsolutiondef} we get 
\begin{equation} \label{szetazero}
S_{\zeta} \Big( \frac{\partial u}{\partial \nu_{+}} - \Lambda_{W,q} f \Big) = 0
\end{equation}
a.e.~in $\Omega'$. This holds also on $\partial \Omega$ by Proposition \ref{layerpotentials} (c), and the uniqueness in the Dirichlet problem for $\Delta$ in $\Omega$ shows that we have \eqref{szetazero} in $\mR^n$. Then \eqref{singlelayerjump} gives that $\frac{\partial u}{\partial \nu_{+}} = \Lambda_{W,q} f$ on $\partial \Omega$.

If $u_1$ and $u_2$ solve (EP) then $u_1|_{\partial \Omega}$ and $u_2|_{\partial \Omega}$ solve (BE), and if $f_1$ and $f_2$ solve (BE) then the corresponding functions defined by \eqref{epsolutiondef} solve (EP). This shows the uniqueness part.
\end{proof}

\begin{remark}
Following Nachman \cite{nachmantwodim}, (BE) is equivalent to 
\begin{equation*}
\left\{ \begin{array}{l} 
(I + S_{\zeta}(\Lambda_{W,q} - \Lambda_{0,0})) f = e^{i\zeta \cdot x} \text{ on } \partial \Omega \\[2pt] 
f \in H^{3/2}(\partial \Omega). 
\end{array} \right.
\end{equation*}
This follows since for $x \in \Omega'$ 
\begin{align*}
D_{\zeta} f(x) &= \int_{\partial \Omega} \frac{\partial G_{\zeta}}{\partial \nu(y)}(x,y) f(y) \,dS(y) = \langle \Lambda_{0,0} G_{\zeta}(x,\,\cdot\,), f \rangle \\
 &= \langle G_{\zeta}(x,\,\cdot\,), \Lambda_{0,0} f \rangle = S_{\zeta} (\Lambda_{0,0} f)(x)
\end{align*}
and letting $x \to \partial \Omega$ nontangentially gives $\frac{1}{2}I + B_{\zeta} = S_{\zeta} \Lambda_{0,0}$.
\end{remark}

Proposition \ref{equivalentproblems} is an immediate consequence of Lemmas \ref{deieequivalence} to \ref{epbeequivalence}. We conclude the section by showing that the operator arising in (BE) is of the form $I + K$ with $K$ compact. This fact and the Fredholm alternative show, for instance, that uniqueness in one of the problems in Proposition \ref{equivalentproblems} implies the existence of a unique solution for all the problems.

\begin{lemma}
Let $\Omega \subseteq \mR^n$, $n \geq 3$, be a bounded domain with $C^{1,1}$ boundary. Then the operator $S_{\zeta} \Lambda_{W,q} - B_{\zeta}- \frac{1}{2} I : H^{3/2}(\partial \Omega) \to H^{3/2}(\partial \Omega)$ is compact.
\end{lemma}
\begin{proof}
Let $f \in H^{3/2}(\partial \Omega)$ and let $u = P_{W,q} f$. If $x \in \Omega$ we have 
\begin{multline*}
-\int_{\Omega} (G_{\zeta}^{\varepsilon}(x,y)\Delta u(y) - u(y)\Delta_y G_{\zeta}^{\varepsilon}(x,y)) \,dy \\
 = \int_{\partial \Omega} \Big( G_{\zeta}^{\varepsilon}(x,y) \frac{\partial u}{\partial \nu}(y) - u(y) \frac{\partial G_{\zeta}^{\varepsilon}(x,y)}{\partial \nu(y)} \Big) \,dS(y).
\end{multline*}
If $\varepsilon \to 0$ we get 
\begin{equation*}
u(x) + \int_{\Omega} G_{\zeta}(x,y) (2W \cdot Du + (W^2 + D \cdot W + q)u) \,dy = (S_{\zeta} \Lambda_{W,q} - D_{\zeta})f(x)
\end{equation*}
a.e.~in $\Omega$. Let then $x \to \partial \Omega$ nontangentially, so that Proposition \ref{layerpotentials} gives 
\begin{equation*}
(S_{\zeta} \Lambda_{W,q} - B_{\zeta} - \frac{1}{2} I)f  = R \int_{\Omega} G_{\zeta}(x,y) (2W \cdot D + (W^2 + D \cdot W + q)) P_{W,q} f(y) \,dy.
\end{equation*}
This reads 
\begin{equation*}
S_{\zeta} \Lambda_{W,q} - B_{\zeta} - \frac{1}{2} I = R G_{\zeta} M J P_{W,q}
\end{equation*}
where $R$ is the trace $H^2(\Omega) \to H^{3/2}(\partial \Omega)$, $G_{\zeta}: L^2(\Omega) \to H^2(\Omega)$ restricts $G_{\zeta} \ast \tilde{u} = e^{i\zeta \cdot x} \Delta_{\zeta}^{-1} e^{-i\zeta \cdot x} \tilde{u}$ to $\Omega$ where $\tilde{u}$ is the extension by zero of $u \in L^2(\Omega)$ to $\mR^n$, $M: H^1(\Omega) \to L^2(\Omega)$ maps $u$ to $2W \cdot Du + (W^2 + D \cdot W + q)u$, and $J$ is the embedding $H^2(\Omega) \to H^1(\Omega)$. All these maps are bounded and $J$ is compact, so the composition is compact.
\end{proof}

\section{Reconstruction of the magnetic field}

The preceding section considered equivalent formulations for problems which give rise to CGO solutions, but did not consider the solvability of any of the problems. The next result, which follows directly from Theorem \ref{normestimates}, shows that if $W$ is continuous then the problems indeed have unique solutions for large $\zeta$.

\begin{prop} \label{equivproblemssolvable}
Assume the conditions in the beginning of Theorem \ref{reconstruction}. Then there exists $C = C(n,\Omega,W,q)$ so that whenever $\abs{\zeta} \geq C$, then each of the problems (DE), (IE), (EP), (BE) has a unique solution.
\end{prop}
\begin{proof}
It is enough to show that (DE) has a unique solution. Now $u = e^{i\zeta \cdot x}(1+\omega)$ solves $H_{W,q} u = 0$ if and only if 
\begin{equation} \label{omegaeq}
(\Delta_{\zeta} + 2W \cdot D_{\zeta} + (W^2 + D \cdot W + q)) \omega = -(2\zeta \cdot W + W^2 + D \cdot W + q).
\end{equation}
From Theorem \ref{normestimates} we know that if $\abs{\zeta} \geq C(n,\Omega,W,q)$ this equation has a unique solution $\omega \in H^1_{\delta}$. Then \eqref{omegaeq} gives that $\Delta_{\zeta} \omega = f$ for some $f \in L^2_{\Omega}$, so that $\omega \in \Delta_{\zeta}^{-1} L^2_{\Omega}$. This shows that $u$ is the unique solution of (DE).
\end{proof}

We will from now on assume the conditions in the beginning of Theorem \ref{reconstruction}. For given $\zeta$ we denote by $u_{\zeta}$ the unique solution of (DE). It follows that if one knows $\Lambda_{W,q}$ then the boundary values $u_{\zeta}|_{\partial \Omega}$ may be reconstructed as the unique solution of the boundary integral equation (BE). The rest of the section will be devoted to showing that the magnetic field $\curl\,W$ may be reconstructed from this knowledge. The first step, similarly as in \cite{nachman}, is to consider a (non-physical) scattering transform.

\begin{definition}
Let $\xi \in \mR^n$ be such that $\abs{\xi}^2$ is not a Dirichlet eigenvalue of $\Delta$ in $\Omega$. Then for any $\zeta \in \mC^n$ which satisfies $\zeta^2 = 0$, $\abs{\zeta} \geq C$, $\re\,\zeta \perp \xi$, and $\im\,\zeta \perp \xi$, we define 
\begin{equation*}
t_{W,q}(\xi,\zeta) = \langle (\Lambda_{W,q} - \Lambda_{0,-\abs{\xi}^2}) (u_{\zeta} |_{\partial \Omega}), e^{-ix \cdot (\xi + \zeta)} |_{\partial \Omega} \rangle.
\end{equation*}
\end{definition}

It is clear from the preceding discussion that $\Lambda_{W,q}$ determines $t_{W,q}$ for the appropriate $\xi, \zeta$. Using the weak formulation of the DN map and the fact that $H_{0,-\abs{\xi}^2} e^{-ix \cdot (\xi+\zeta)} = 0$ in $\Omega$, one easily sees that 
\begin{equation} \label{twqeq}
t_{W,q}(\xi,\zeta) = \int_{\Omega} e^{-ix \cdot \xi} (2 (\zeta \cdot W) u_0 + W \cdot Du_0 + (W^2 + \xi \cdot W + \abs{\xi}^2 + q)u_0) \,dx,
\end{equation}
where we write $u_0 = e^{-i\zeta \cdot x} u_{\zeta}$.

We know from (DE) that $u_0 = 1 + \omega$ with $\omega \in \Delta_{\zeta}^{-1} L^2_{\Omega}$. If one had $\norm{\omega}_{L^2(\Omega)} \to 0$ as $\abs{\zeta} \to \infty$ one could divide \eqref{twqeq} by $\abs{\zeta}$ and let $\abs{\zeta} \to \infty$, which would then give essentially the Fourier transform of $\curl\,W$. However, $\omega$ is obtained by solving \eqref{omegaeq} where the $L^2_{\delta+1}$ norm of the right hand side is $O(\abs{\zeta})$ instead of $o(\abs{\zeta})$, so one gets that $\norm{\omega}_{L^2(\Omega)}$ is bounded but may not be small when $\abs{\zeta}$ is large (in fact Lemma \ref{modifiedcgolemma} shows that $\omega \to e^{i\phi} - 1$ in $L^2(\Omega)$ as $\abs{\zeta} \to \infty$, where $\phi$ is defined below).

To deal with this difficulty we write the solution $u_{\zeta}$ in a different form where one gets a small remainder term for $\abs{\zeta}$ large. For this we employ a decomposition 
\begin{equation*}
W = W^{\sharp} + W^{\flat}
\end{equation*}
where $W^{\sharp} = W \ast \varphi_{\varepsilon}$ with $\varphi_{\varepsilon}$ the usual mollifier, and we make the choice 
\begin{equation*}
\varepsilon = \abs{\zeta}^{-\sigma}
\end{equation*}
with $0 < \sigma < 1/2$. Then $W^{\sharp}$ is $C^{\infty}$ and 
\begin{align}
\norm{W^{\sharp}}_{W^{1,\infty}} &= o(\abs{\zeta}^{\sigma}), \\
\norm{W^{\sharp}}_{W^{2,\infty}} &= o(\abs{\zeta}^{2\sigma}), \\
\norm{W^{\flat}}_{L^{\infty}} &= o(1)
\end{align}
as $\abs{\zeta} \to \infty$.

We also write $\zeta = s \mu$ where $\mu = \gamma_1 + i \gamma_2$, $\abs{\gamma_j} = 1$, $\gamma_1 \perp \gamma_2$. Finally, we fix $\chi \in C^{\infty}_c(\mR^n)$ with $\chi = 1$ in $B(0,M/2)$, $\chi = 0$ outside of $B(0,M)$, and $\closure{\Omega} \subseteq B(0,M/2)$.

\begin{lemma} \label{modifiedcgolemma}
Fix $\theta > 0$ with $\sigma + \theta < 1/2$. For $\abs{\zeta}$ large enough, the CGO solution $u_{\zeta}$ of $H_{W,q} u_{\zeta} = 0$ in $\mR^n$ may be written in the form 
\begin{equation} \label{modifiedcgo}
u_{\zeta} = e^{i\zeta \cdot x}(\omega_0 + \omega)
\end{equation}
where $\omega_0 = e^{i\chi_{\zeta} \phi^{\sharp}}$ with 
\begin{align} \label{phisharpdef}
\chi_{\zeta}(x) &= \chi(x/\abs{\zeta}^{\theta}), \\
\phi^{\sharp}(x) &= N_{\mu}^{-1}(-\mu \cdot W^{\sharp}),
\end{align}
and $\omega \in H^1_{\delta}$ with $\norm{\omega}_{L^2_{\delta}} = o(1), \norm{\omega}_{H^1_{\delta}} = o(\abs{\zeta})$ as $\abs{\zeta} \to \infty$.
\end{lemma}
\begin{proof}
We first show that the equation $H_{W,q} u = 0$ in $\mR^n$ has a solution of the form \eqref{modifiedcgo} with the required properties. This will be the case if $\omega$ satisfies 
\begin{equation} \label{modifiedomegaeq}
(\Delta_{\zeta} + 2W \cdot D_{\zeta} + G) \omega = -f
\end{equation}
where $G = W^2 + D \cdot W + q \in L^{\infty}_{\Omega}$ and 
\begin{multline*}
f = (\Delta_{\zeta} + 2W \cdot D_{\zeta} + G) \omega_0 = e^{i\chi_{\zeta}\phi^{\sharp}} \Big[ i\chi_{\zeta}\Delta\phi^{\sharp} + 2i D\chi_{\zeta} \cdot D\phi^{\sharp} + i\phi^{\sharp}\Delta\chi_{\zeta} \\
 + (\chi_{\zeta} \nabla\phi^{\sharp} + \phi^{\sharp}\nabla \chi_{\zeta})^2 + 2\zeta \cdot (\nabla \chi_{\zeta}) \phi^{\sharp} + 2\zeta \cdot (\nabla \phi^{\sharp}) \chi_{\zeta} \\
  + 2W \cdot (\nabla \chi_{\zeta})\phi^{\sharp} + 2W \cdot (\nabla \phi^{\sharp}) \chi_{\zeta} + 2W^{\sharp} \cdot \zeta + 2W^{\flat} \cdot \zeta + G \Big].
\end{multline*}
We need to know the behaviour of $\norm{f}_{L^2_{\delta+1}}$ as $\abs{\zeta}$ grows. The choice of $\phi^{\sharp}$ implies 
\begin{equation*}
2\zeta \cdot \nabla \phi^{\sharp} + 2W^{\sharp} \cdot \zeta = 0.
\end{equation*}
Since $W^{\sharp} = \chi_{\zeta} W^{\sharp}$ this removes the worst two terms from $f$, and one obtains in terms of $L^2_{\delta+1}$ norms 
\begin{multline} \label{fbigestimate}
\norm{f} \leq C \Big[ \norm{\chi_{\zeta}\Delta\phi^{\sharp}} + \norm{\nabla \chi_{\zeta} \cdot \nabla \phi^{\sharp}} + \norm{\phi^{\sharp}\Delta\chi_{\zeta}}
 + \norm{\abs{\chi_{\zeta} \nabla\phi^{\sharp}}^2} + \norm{\abs{\phi^{\sharp}\nabla \chi_{\zeta}}^2} \\
 + \abs{\zeta}^{1-\theta} \norm{(\nabla \chi(x/\abs{\zeta}^{\theta})) \phi^{\sharp}} + \norm{W \cdot (\nabla \chi_{\zeta})\phi^{\sharp}} + \norm{W \cdot (\nabla \phi^{\sharp}) \chi_{\zeta}} \\
 + \norm{W^{\flat} \cdot \zeta} + \norm{G} \Big].
\end{multline}
Lemma \ref{nmuinv1} implies 
\begin{equation} \label{nmuinvphisharpest}
\abs{\partial^{\alpha} \phi^{\sharp}(x)} \leq C \abs{\zeta}^{\sigma\abs{\alpha}} \br{x_T}^{-1} \chi_{B(0,M)}(x_{\perp})
\end{equation}
where $x_T$ is the projection of $x$ to $\mathrm{span}\{\gamma_1,\gamma_2\}$ and $x_{\perp} = x-x_T$. Then for instance 
\begin{align*}
\norm{\chi_{\zeta} \Delta \phi^{\sharp}} &= \Big( \int_{\mR^n} \br{x}^{2(\delta+1)} \chi_{\zeta}(x)^2 \abs{\Delta \phi^{\sharp}(x)}^2 \,dx \Big)^{1/2} \\
 &\leq C \abs{\zeta}^{2\sigma} \Big( \int_{\abs{x_T} \leq M \abs{\zeta}^{\theta}, \abs{x_{\perp}} \leq M} \br{x}^{2(\delta+1)} \br{x_T}^{-2} \,dx \Big)^{1/2} \\
 &\leq C \abs{\zeta}^{2\sigma} \Big( \int_{\abs{x_T} \leq M \abs{\zeta}^{\theta}} \br{x_T}^{2\delta} \,dx_T \Big)^{1/2} \leq C \abs{\zeta}^{2\sigma+(\delta+1)\theta}.
\end{align*}
This has the worst behaviour of the first five terms of \eqref{fbigestimate} since derivatives hitting $\chi_{\zeta}$ bring decay in $\abs{\zeta}$ and the other terms involve only first derivatives of $\phi^{\sharp}$. A similar computation shows that the sixth term is $O(\abs{\zeta}^{1-\theta+(\delta+1)\theta}) = O(\abs{\zeta}^{1+\delta\theta})$. One has $W \cdot \nabla \chi_{\zeta} = 0$ and $W \chi_{\zeta} = W$ for large $\abs{\zeta}$ so the seventh and eight terms are $0$ and $O(\abs{\zeta}^{\sigma})$, respectively. The final two terms are $o(\abs{\zeta})$ and $O(1)$, respectively, since $\norm{W^{\flat}}_{L^{\infty}} \to 0$ as $\abs{\zeta} \to \infty$. Using the choices of $\sigma$ and $\theta$ and the fact that $-1 < \delta < 0$, we obtain $\norm{f}_{L^2_{\delta+1}} = o(\abs{\zeta})$. The solution $\omega$ of \eqref{modifiedomegaeq} has the desired properties by Theorem \ref{normestimates}.

It remains to show that $u$ given by \eqref{modifiedcgo} is the CGO solution. One has $u = e^{i\zeta \cdot x}(1+\tilde{\omega})$ where 
\begin{equation*}
\tilde{\omega} = e^{i\chi_{\zeta} \phi^{\sharp}} - 1 + \omega.
\end{equation*}
Now $e^{i\chi_{\zeta} \phi^{\sharp}} - 1$ is in $H^1_{\delta}$ for any $\delta < 0$, since for instance $e^{i\chi_{\zeta} \phi^{\sharp}} - 1 = O(\abs{\chi_{\zeta} \phi^{\sharp}})$ and \eqref{nmuinvphisharpest} implies 
\begin{multline*}
\norm{e^{i\chi_{\zeta} \phi^{\sharp}} - 1}_{L^2_{\delta}} \leq C \Big( \int_{\mR^n} \br{x}^{2\delta} \chi_{\zeta}(x)^2 \abs{\phi^{\sharp}(x)}^2 \,dx \Big)^{1/2} \\
 \leq C \Big( \int_{\mR^2} \br{x_T}^{2\delta-2} \,dx_T \Big)^{1/2} < \infty.
\end{multline*}
Also $\omega \in H^1_{\delta}$ so $\tilde{\omega} \in H^1_{\delta}$ for $-1 < \delta < 0$. This and $H_{W,q} u = 0$ imply $\tilde{\omega} \in \Delta_{\zeta}^{-1} L^2_{\Omega}$, so $u$ is indeed the unique solution of (DE) given by Proposition \ref{equivproblemssolvable}.
\end{proof}

We may now plug in $u_0 = \omega_0 + \omega$ from \eqref{modifiedcgo} to \eqref{twqeq}. The estimates for $\omega$ and the form of $\omega_0$ imply that 
\begin{equation*}
R_{W,q}(\xi,\mu) = \lim_{s \to \infty} s^{-1} t_{W,q}(\xi,s\mu) = 2 \int e^{-ix \cdot \xi}e^{i\phi}(\mu \cdot W) \,dx
\end{equation*}
where $\phi = N_{\mu}^{-1}(-\mu \cdot W)$. This shows that we may recover a nonlinear Fourier transform $R_{W,q}(\xi,\mu)$ of $\mu \cdot W$ from the knowledge of $\Lambda_{W,q}$, for any $\mu$ and $\xi$ with $\xi \cdot \mu = 0$ and $\abs{\xi}^2$ not a Dirichlet eigenvalue of $\Delta$ in $\Omega$.

The next argument, due to Eskin and Ralston \cite{eskinralston}, shows that this nonlinear Fourier transform is in fact just an ordinary Fourier transform. Similar ideas appear in Sun \cite{sun}.

\begin{lemma}
One has 
\begin{equation*}
R_{W,q}(\xi,\mu) = 2 \int e^{-ix \cdot \xi} (\mu \cdot W) \,dx.
\end{equation*}
\end{lemma}
\begin{proof}
It is enough to prove this for $\mu = e_1 + i e_2$, so that $\xi = (0,\xi'')$ and $(\partial_1 + i \partial_2) \phi = -(W_1 + i W_2)$. Then 
\begin{align*}
R_{W,q}(\xi,\mu) &= 2 \int_{\mR^n} e^{-ix'' \cdot \xi''} e^{i\phi} (-(\partial_1 + i \partial_2) \phi) \,dx \\
 &= 2 \int_{\mR^{n-2}} e^{-ix'' \cdot \xi''} h(x'') \,dx''
\end{align*}
where 
\begin{align*}
h(x'') &= i \int_{\mR^2} (\partial_1 + i \partial_2)(e^{i\phi(x',x'')}) \,dx' = \lim_{R \to \infty} i \int_{\abs{x'} \leq R} (\partial_1 + i \partial_2)(e^{i\phi(x',x'')}) \,dx' \\
 &= \lim_{R \to \infty} i \int_{\abs{x'} = R} e^{i\phi(x',x'')} (\nu_1 + i \nu_2) \,dS(x').
\end{align*}
The integrals are well defined by standard approximation arguments. Now $e^{i\phi} = 1 + i\phi + O(\abs{i\phi}^2) = 1 + i\phi + O(\abs{x'}^{-2})$ by Lemma \ref{nmuinv1}, and 
\begin{eqnarray*}
 & \int_{\abs{x'} = R} (\nu_1 + i \nu_2) \,dS(x') = \int_{\abs{x'} \leq R} (\partial_1 + i \partial_2)(1) \,dx' = 0, & \\
 & \Big\lvert \int_{\abs{x'} = R} O(\abs{x'}^{-2})(\nu_1 + i \nu_2) \,dS(x') \Big\rvert \leq \frac{C}{R} \to 0 \text{ as } R \to \infty, & 
\end{eqnarray*}
so we have 
\begin{align*}
h(x'') &= -\lim_{R \to \infty} \int_{\abs{x'} = R} \phi(x',x'')(\nu_1 + i \nu_2) \,dS(x') \\
 &= -\lim_{R \to \infty} \int_{\abs{x'} \leq R} (\partial_1 + i \partial_2) \phi \,dx' = \int_{\mR^2} (e_1 + i e_2) \cdot W \,dx'.
\end{align*}
This gives the claim.
\end{proof}

We now show that $R_{W,q}(\xi,\mu)$ determines $\curl\,W$, or $D_j W_k - D_k W_j$ for any $j \neq k$. Let $\xi \in \mR^n$ be such that $\abs{\xi}^2$ is not a Dirichlet eigenvalue of $\Delta$ in $\Omega$. If one of $\xi_j, \xi_k$ is nonzero choose $\gamma$ to be the unit vector with direction $\xi_j e_k - \xi_k e_j$, so that $\xi \cdot \gamma = 0$. Since $n \geq 3$ we may choose a unit vector $\tilde{\gamma}$ with $\gamma \cdot \tilde{\gamma} = \xi \cdot \tilde{\gamma} = 0$. Letting $\mu = \gamma + i \tilde{\gamma}$ we see that $\Lambda_{W,q}$ determines $R_{W,q}(\xi,\mu) + R_{W,q}(\xi,\overline{\mu})$, which determines 
\begin{equation*}
\int e^{-ix \cdot \xi} (\xi_j W_k - \xi_k W_j) \,dx = (D_j W_k - D_k W_j)\ehat(\xi)
\end{equation*}
in fact for any $\xi \in \mR^n$ such that $\abs{\xi}^2$ is not an eigenvalue. Since there are countably many eigenvalues and since $D_j W_k - D_k W_j$ is compactly supported, so the Fourier transform is analytic, we recover $D_j W_k - D_k W_j$.

\section{Reconstruction of the electric potential}

Finally, we make the additional assumptions that $W$ is $C^{1+\varepsilon}$ and $\partial \Omega$ is $C^{2+\varepsilon}$ for some $\varepsilon > 0$, and we indicate how to recover $q$ from $\Lambda_{W,q}$. From the preceding section, we may assume that the magnetic field $\curl\,W$ is known. The next step is to construct a certain magnetic potential with this magnetic field.

\begin{lemma}
One can construct $\tilde{W} \in C^{1+\varepsilon}(\Omega ; \mC^n)$ with $\curl\,\tilde{W} = \curl\,W$ and $\tilde{W}|_{\partial \Omega} = 0$.
\end{lemma}
\begin{proof}
For the following concepts we refer to \cite{schwarz}. We write $X \Lambda^k(\Omega)$ for a differential $k$-form in $\Omega$ with coefficient functions in $X$, and $t\eta$ and $n\eta$ for the tangential and normal traces on $\partial \Omega$, respectively, of a form $\eta$. Let $\omega = \sum_{j=1}^n W_j dx_j \in C^{1+\varepsilon} \Lambda^1(\Omega)$ be the $1$-form corresponding to $W$, and let $\chi = d\omega$. We start with solving the boundary value problem 
\begin{equation} \label{exteriordifferentialbvp}
\left\{ \begin{array}{rll}
d\tilde{\omega} &\!\!\!= \chi & \quad \text{in } \Omega, \\
t\tilde{\omega} &\!\!\!= 0 & \quad \text{on } \partial \Omega
\end{array} \right.
\end{equation}
for $\tilde{\omega}$ in $C^{1+\varepsilon}\Lambda^1(\Omega)$. In the case of smooth domains and $L^p$ Sobolev spaces this problem is considered in \cite{schwarz}, where the solution is reduced to the Hodge decomposition of $\chi$. For a general form $\chi \in L^2\Lambda^k(\Omega)$ this decomposition reads 
\begin{equation} \label{hodgedecomposition}
\chi = d\alpha + \delta \beta + \kappa
\end{equation}
where $\alpha \in H^1 \Lambda^{k-1}(\Omega)$ with $t\alpha = 0$, $\beta \in H^1 \Lambda^{k+1}(\Omega)$ with $n\beta = 0$, and $\kappa$ is a harmonic field meaning that $d\kappa = \delta \kappa = 0$. Here $\delta$ is the codifferential. Further, the three summands in  \eqref{hodgedecomposition} are uniquely determined and mutually orthogonal with respect to the natural $L^2$ inner product.

The specific form of $\chi$ above implies that the decomposition \eqref{hodgedecomposition} reduces to $\chi = d\alpha$, where one may choose $\alpha = \delta \phi$ where $\phi$ is the Dirichlet potential of $\chi$ (see \cite{schwarz}, Section 2.2). Under the present assumptions of $C^{2+\varepsilon}$ boundary and $C^{\varepsilon}$ regularity of $\chi$, Theorem 7.7.4 in \cite{morrey} implies that $\delta \phi \in C^{1+\varepsilon}$. It is then easy to check that $\tilde{\omega} = \delta \phi$ solves \eqref{exteriordifferentialbvp}. We note that the Dirichlet potential $\phi$ is obtained constructively using the explicit integral formula for the corresponding Green operator, as in \cite{mitreamitrea}.

Letting $\tilde{W}$ be the vector field corresponding to $\tilde{\omega}$, we have that $\tilde{W} \in C^{1+\varepsilon}(\Omega ; \mC^n)$, $\curl\,\tilde{W} = \curl\,W$, and the tangential components of $\tilde{W}$ vanish on $\partial \Omega$. We may further replace $\tilde{W}$ by $\tilde{W} + \nabla p$ where $p \in C^{2+\varepsilon}(\Omega)$ satisfies $p|_{\partial \Omega} = 0$ and $\frac{\partial p}{\partial \nu}|_{\partial \Omega} = -\tilde{W} \cdot \nu$, and $p$ is constructed similarly as in Lemma 5.8 of \cite{salo}. This completes the proof.
\end{proof}

With $\tilde{W}$ as above, we conclude that $\tilde{W} = W + \nabla p$ where $p \in W^{1,\infty}(\Omega)$ and $p|_{\partial \Omega} = 0$. An easy argument for this under the present regularity assumptions is obtained by extending $\tilde{W} - W$ by zero to $\mR^n$ as a Lipschitz vector field, and by noting that $\tilde{W} - W = \nabla p$ in a large ball with $p \in W^{2,\infty}$ and $\nabla p|_{\partial \Omega} = 0$. The assumption on the topology of $\Omega$ ensures that $\partial \Omega$ is connected, so that $p$ is constant on $\partial \Omega$, and one may substract the constant to get $p|_{\partial \Omega} = 0$. Gauge equivalence then implies that $\Lambda_{\tilde{W},0} = \Lambda_{W,0}$.

Fix $\xi \in \mR^n \smallsetminus \{0\}$ and take unit vectors $\gamma_j$ such that $\{\xi,\gamma_1,\gamma_2\}$ form an orthogonal set. For $s > 0$ define complex vectors 
\begin{align*}
\zeta_1 &= -\frac{\xi}{2} + s\sqrt{1-\frac{\abs{\xi}^2}{4s^2}} \gamma_1 + is \gamma_2, \\
\zeta_2 &= -\frac{\xi}{2} - s\sqrt{1-\frac{\abs{\xi}^2}{4s^2}} \gamma_1 - is \gamma_2.
\end{align*}
Using the notation of Lemma \ref{modifiedcgolemma}, the equations $H_{W,q} u = 0$ and $H_{-W,0} v = 0$ have unique CGO solutions $u = u_{\zeta_1}$ and $v = v_{\zeta_2}$, which have the form 
\begin{align*}
u_{\zeta_1} &= e^{i\zeta_1 \cdot x}(e^{i\chi_{\zeta_1} \phi^{\sharp}} + \omega_1), \\
v_{\zeta_2} &= e^{i\zeta_2 \cdot x}(e^{-i\chi_{\zeta_2} \phi^{\sharp}} + \omega_2)
\end{align*}
where $\norm{\omega_j}_{L^2_{\delta}} \to 0$ as $s \to \infty$.

We define a new scattering transform 
\begin{equation*}
\tilde{t}(\xi) = \langle (\Lambda_{W,q} - \Lambda_{W,0}) (u_{\zeta_1}|_{\partial \Omega}), v_{\zeta_2} |_{\partial \Omega} \rangle.
\end{equation*}
Since $\Lambda_{W,q}$ and $\Lambda_{-W,0}$ are known, one may construct the boundary values of $u_{\zeta_1}$ and $v_{\zeta_2}$ as solutions of boundary integral equations as in Section 5, and thus $\tilde{t}$ is known. The definition of DN maps implies 
\begin{equation*}
\tilde{t}(\xi) = \int q u_{\zeta_1} v_{\zeta_2} \,dx,
\end{equation*}
and so $\tilde{t}(\xi) \to \int e^{-ix \cdot \xi} q(x) \,dx$ as $s \to \infty$. This is the Fourier transform of $q$, and we have recovered the electric potential.

\subsection*{Acknowledgements}

I would like to thank Anders Melin, James Ralston, and Michael Taylor for useful suggestions. Part of the research was done while visiting the University of Washington, and I would like to thank Gunther Uhlmann for his generous support and helpful discussions. Financial support from the Finnish Academy of Science and Letters, Vilho, Yrj\"o and Kalle V\"ais\"al\"a Foundation is gratefully acknowledged.

%\nocite{*}
\addcontentsline{toc}{chapter}{Bibliography}
\bibliography{magnetic}
\bibliographystyle{hamsplain}

\end{document}